\documentclass[12pt,reqno]{amsart}
\usepackage{amsmath,amsfonts,amsthm,amsopn,amssymb}
\usepackage{cite,marginnote}
\usepackage{color,enumitem,graphicx}
\usepackage[colorlinks=true,urlcolor=blue,
citecolor=red,linkcolor=blue,linktocpage,pdfpagelabels,
bookmarksnumbered,bookmarksopen]{hyperref}
\usepackage[english]{babel}

\usepackage[left=2.9cm,right=2.9cm,top=2.8cm,bottom=2.8cm]{geometry}

\numberwithin{equation}{section}

\makeindex

\newtheorem{Remark}{Remark}[section]

\newtheorem{lemma}{Lemma}[section]
\newtheorem{iteration lemma}{iteration Lemma}[section]

\newcommand{\bt}{\begin{theorem}}
\newcommand{\et}{\end{theorem}}
\newcommand{\bl}{\begin{lemma}}
\newcommand{\el}{\end{lemma}}
\newcommand{\bd}{\begin{definition}}
\newcommand{\ed}{\end{definition}}
\newcommand{\bc}{\begin{corollary}}
\newcommand{\ec}{\end{corollary}}
\newcommand{\bp}{\begin{proof}}
\newcommand{\ep}{\end{proof}}
\newcommand{\bx}{\begin{example}}
\newcommand{\ex}{\end{example}}
\newcommand{\bi}{\begin{exercise}}
\newcommand{\ei}{\end{exercise}}
\newcommand{\bo}{\begin{proposition}}
\newcommand{\eo}{\end{proposition}}
\newcommand{\br}{\begin{remark}}
\newcommand{\er}{\end{remark}}
\newcommand{\beq}{\begin{equation}}
\newcommand{\eeq}{\end{equation}}
\newcommand{\ba}{\begin{align}}
\newcommand{\ea}{\end{align}}
\newcommand{\bn}{\begin{enumerate}}
\newcommand{\en}{\end{enumerate}}
\newcommand{\bg}{\begin{align*}}
\newcommand{\bcs}{\begin{cases}}
\newcommand{\ecs}{\end{cases}}

\newcommand{\bean}{\begin{eqnarray*}}
\newcommand{\eean}{\end{eqnarray*}}

\def\bd{\mathrm{bd}\,}

\title[Nodal solutions for quasilinear Schr\"{o}dinger equations]{Nodal solutions for quasilinear Schr\"{o}dinger equations with asymptotically 3-linear nonlinearity}

\author[H. Zhang]{Hui Zhang}
\author[F.~J.~Meng]{Fengjuan Meng}
\author[J. J. Zhang]{Jianjun Zhang}

\address[H.\ Zhang]{\newline\indent Department of  Mathematics,
Jinling Institute of Technology,
\newline\indent
Nanjing 211169, PR China
\newline\indent and
\newline\indent Department of Mathematics,
Nanjing University,
\newline\indent
Nanjing 210093, PR China}
\email{\href{mailto:huihz0517@126.com}{huihz0517@126.com}}

\address[F. J.\ Meng]{\newline\indent School of Mathematics and Physics
\newline\indent
Jiangsu University of Technology
\newline\indent
Changzhou 213001, China}
\email{\href{mailto:fjmeng@jsut.edu.cn}{fjmeng@jsut.edu.cn}}

\address[J. J. \ Zhang]{\newline\indent College of Mathematics and Statistics, Chongqing Jiaotong University,
\newline\indent
Chongqing 400074, PR China}
\email{\href{mailto:zhangjianjun09@tsinghua.org.cn}{zhangjianjun09@tsinghua.org.cn}}

\thanks{(1) Corresponding author: Jianjun Zhang ({\tt zhangjianjun09@tsinghua.org.cn})}

\thanks{(2) Hui Zhang was supported by China Postdoctoral Science Foundation (No. 2021M691527). Fengjuan Meng was supported by QingLan Project of Jiangsu Province. Jianjun Zhang was supported by National Natural Science Foundation of China (No. 11871123).}

\subjclass[2000]{35J20; 35J60; 35Q55}
\keywords{Quasilinear Schr\"{o}dinger equation; Asymptotically 3-linear nonlinearity; Nodal solution; Variational method.}

\begin{document}

\begin{abstract}
In this paper, we are concerned with the quasilinear Schr\"{o}dinger equation
\begin{equation*}
-\Delta u+V(x)u-u\Delta(u^2)=g(u),\ \
x\in \mathbb{R}^{N},
\end{equation*}
where $N\geq3$, $V$ is radially symmetric and nonnegative, and $g$ is asymptotically 3-linear at infinity. In the case of $\inf_{\mathbb{R}^N}V>0$, we show the existence of a least energy sign-changing solution with exactly one node, and for any integer $k>0$, there are a pair of sign-changing solutions with $k$ nodes. Moreover, in the case of $\inf_{\mathbb{R}^N}V=0$, the problem above admits a least energy sign-changing solution with exactly one node. The proof is based on variational methods. In particular, some new tricks and the method of sign-changing Nehari manifold depending on a suitable restricted set are introduced to overcome the difficulty resulting from the appearance of asymptotically 3-linear nonlinearities.
\end{abstract}

\maketitle

\section{Introduction and main results}
\setcounter{equation}{0}
Quasilinear Schr\"{o}dinger equations of the form
\begin{equation}\label{1.1}
i\partial_tz=-\triangle z+W(x)z-m(|z|^2)z-\kappa\triangle
\phi(|z|^2)\phi'(|z|^2)z, \hskip 3 cm
\end{equation}
have been derived as models of several physical phenomena and extensively studied in recent years, where
$N\geq3$, $z:\mathbb{R}\times \mathbb{R}^N\rightarrow \mathbb{C},$
$W:\mathbb{R}^N\rightarrow \mathbb{R}$ is a given potential,
$m,\phi:\mathbb{R}^+\rightarrow \mathbb{R}$ are suitable functions, and
$\kappa$ is a real constant. For example, if $\kappa=0,$ (\ref{1.1}) turns out to be  semilinear Schr\"{o}dinger equations, which have been widely
investigated, we refer the readers to \cite{BERE,Rap,WM1}. In the case
$\phi(s)=s$, as a model of the time evolution of the condensate wave
function in super-fluid film, has been studied by Kurihara in
\cite{Ku}. While for the case $\phi(s)=\sqrt{1+s},$ the equations are the models
of the self-channeling of a high-power ultra short laser in matter, see \cite{ChenS}. For more physical
applications, we refer to \cite{BL,LPT,PS} and references therein.
 Here we are interested in the case $\phi(s) = s$ and $\kappa = 1$.
Looking for standing wave solutions of (\ref{1.1}), that is,
solutions of the form $z(t, x) =e^{-iE t}u(x)$, $E\in\mathbb{R},$ we
get the following equation
\begin{equation}\label{1.2}
-\triangle u+V(x)u-u \Delta (u^2)={g}(u), \ \  x\in \mathbb{R}^N,\hskip 3 cm
\end{equation} where $N\geq3$, $V(x)=W(x)-E$
 and ${g}(u)=m(u^2)u$.

There are a few results on the existence and multiplicity of sign-changing solutions of (\ref{1.2}). In \cite{LWW}, Liu et. al.
considered the equation
\begin{equation}\label{1.3}
-\triangle u+V(x)u-u \Delta (u^2)=|u|^{p-2}u, \ \  x\in \mathbb{R}^N,
\end{equation}
where $p\in[4,22^*)$ with $2^*=\frac{2N}{N-2}$, and $V(x)\in C(\mathbb{R}^N)$ satisfies
\vskip 0.2 true cm
\noindent(V$'_1$)\ $0<\inf_{\mathbb{R}^N}V(x)\leq \lim_{|x|\rightarrow+\infty}V(x):=V_\infty$,\ \text{and} \ $V(x)\leq V_\infty-\frac{A}{1+|x|^m}$ \text{for} \ $|x|\geq M$,
\vskip 0.2 true cm
\noindent where $A, M, m>0$. They proved the existence of a least energy sign-changing solution for (\ref{1.3}) by using an approximating sequence of problems on a Nehari manifold defined in an appropriate subset of $H^1(\mathbb{R}^N)$. In \cite{DPW1}, Deng et. al. treated the equation (\ref{1.3}) with $p\in(4,22^*)$ and showed that, for any given integer $k\geq0$, there is a pair of
sign-changing solutions with $k$ nodes by using a minimization argument and  an energy comparison method. Later, Deng et. al. \cite{DPW2} extended the results in \cite{DPW1} to the critical growth case. Recently, Yang et. al. \cite{YANG1} dealt with the critical or supercritical growth case
$$-\triangle u+V(x)u-u \Delta (u^2)=a(x)[g(u)+|u|^{p-2}u], \ \  x\in \mathbb{R}^N,$$
where $p\geq 22^*$.  By assuming that $a(x)>0$ a.e. in $\mathbb{R}^N$ and $V$ satisfies one of the following conditions
\vskip 0.2 true cm
\noindent(V$'_2$) $V(x)\geq V_0>0$ for all $x\in\mathbb{R}^N$,\ $V(x)=V(|x|)$ and $V\in L^\infty(\mathbb{R}^N)$,
\vskip 0.1 true cm
\noindent(V$'_3$) $V(x)\geq V_0>0$ for all $x\in\mathbb{R}^N$,\ $\lim_{|x|\rightarrow\infty}V(x)=+\infty$,
\vskip 0.2 true cm
\noindent Yang et. al.  \cite{YANG1} obtained a least energy sign-changing solution and the proof was based on the method of Nehari manifold, deformation arguments and
$L^\infty$-estimates. Later, Liu et al. \cite{LIULIUWANG} considered  a parameter-dependent quasilinear equation
\begin{equation*}\aligned
\left\{ \begin{array}{lll}
\Delta u+\frac12u\Delta (u^2)+\lambda|u|^{r-2}u=0,\ \quad\ & \text{in}\quad \Omega,\\
u=0, \quad  \ &\text{on}\quad\partial\Omega,
\end{array}\right.\endaligned
\end{equation*}
where $r\in (2,4)$. By means of a perturbation approach, they proved the existence of more and more
sign-changing solutions with both positive and negative energies when $\lambda$ becomes large. Recently, without any parameter-dependence, Zhang et al. \cite{ZLTZ} showed the existence and multiplicity of sign-changing solutions of (\ref{1.2}) with sub-cubic or cubic nonlinearity by virtue of a new perturbation approach and the method of invariant sets
of descending flow.

Observe that no information about nodal domains is available for sign-changing solutions obtained in \cite{LIULIUWANG,ZLTZ} for quasilinear Schr\"{o}dinger equation with sub-cubic or cubic nonlinearity. In this paper, we are interested in sign-changing solutions with nodal domains. As far as we know, little has done on the existence and multiplicity of nodal solutions of (\ref{1.2}) with asymptotically 3-linear nonlinearity.  In addition, in all the works mentioned above, the potential was assumed to be non-vanishing, i.e. $\inf_{\mathbb{R}^N}V>0$. So it is quite natural to ask whether there is still one or more nodal solutions if $g$ is asymptotically 3-linear, and whether nodal solutions still exist if $V$ is vanishing, i.e. $\inf_{\mathbb{R}^N}V=0$. We shall give affirmative answers.

On one hand, we consider the case that $V$ is non-vanishing and  assume
  \vskip 0.2 true cm
 \noindent(V$_1$) $V\in C(\mathbb{R}^{N})\cap L^\infty(\mathbb{R}^{N})$, $V(x)=V(|x|)$ and $V(x)\geq0$ on $\mathbb{R}^{N}$.
 \vskip 0.2 true cm
 \noindent(V$_2$) $\inf_{x\in\mathbb{R}^N}V(x)>0$.
\vskip 0.2 true cm \noindent($g_1$) $g\in C(\mathbb{R})$ and $g(t)=o(t^3)$\ as\ $t\rightarrow0$.
\vskip 0.2 true cm \noindent{($g_2$)} $\lim_{|t|\rightarrow+\infty}\frac{g(t)}{t^3}=1$.
\vskip 0.2 true cm \noindent{($g_3$)}\ $t\mapsto \frac{g(t)}{t^3}$ is decreasing on $(-\infty,0)$ and increasing on $(0,+\infty)$.
 \vskip 0.2 true cm \noindent{($g_4$)}\ $\lim_{|t|\rightarrow+\infty}\bigl[\frac14 g(t)t-G(t)\bigr]=+\infty$, where $G(t):=\int^t_0g(s)ds$.
\begin{Remark} An example of $g$ is $g(t)=\frac{t^5}{1+t^2}$. By direct calculations, we have $G(t)=\frac14t^4-\frac12t^2+\frac12\ln(1+t^2)$. It is easy to see that $g$ satisfies the assumptions (g$_1$)-(g$_4$).
\end{Remark}
\vskip 0.1 true cm

\noindent\textbf{Theorem 1.1.} {\it Let (V$_1$), (V$_2$) and (g$_1$)-(g$_4$) hold. Then the equation (\ref{1.2}) admits a least
energy sign-changing solution, which changes sign exactly once.}
\vskip 0.1 true cm
\noindent\textbf{Theorem 1.2.} {\it Let (V$_1$), (V$_2$) and (g$_1$)-(g$_4$) hold. Then for any integer $k>0$, there
exists a pair $u^{\pm}_k$
of radial solutions of (\ref{1.2}) with the following properties:\\
(1) $u^{-}_k(0)<0<u^{+}_k(0)$;\\
(2) $u^{\pm}_k$ have exactly $k$ nodes $0<\rho^{\pm}_1<...<\rho^{\pm}_k<+\infty$ and $u^{\pm}_k(\rho^{\pm}_k)=0$.}

\begin{Remark} We give some remarks for Theorem 1.2.\\
(1) The functional energy of $u^{\pm}_k$
are at least $(k+1)l$, where $l$ is the least energy corresponding to the equation (\ref{1.2}), see Remark \ref{r4.2} for details. So $u^{\pm}_k$
are high energy sign-changing
solutions if $k\geq1$.\\
(2) We can not claim that $u^+_k=-u^-_k$
since the nodes $\rho^+_i$ may not equal to $\rho^-_i$
for $1\leq i\leq k$.
\end{Remark}

On the other hand, we consider the case that $V$ is vanishing. For any
rotationally symmetric domain $\Omega\subset\mathbb{R}^N$, we denote $H^1_{0,rad}(\Omega)$ as the Hilbert space of radial functions in $H^1_0(\Omega)$. We
suppose\\
\noindent(V$_3$) There exists $a>0$ such that $\mathcal{V}_a:=\{x\in\mathbb{R}^N: V(x)<a\}$ is nonempty and has finite measure.

\noindent(V$_4$) $\Omega_0=\{x\in\mathbb{R}^N: V(x)=0\}\neq\emptyset$ and there exist two bounded, smooth and rotationally
symmetric domains $\Omega_1,\ \Omega_2\subset\Omega_0$ such that $\overline{\Omega_1}\cap\overline{\Omega_2}=\emptyset$.

\noindent($g_5$) $g\in C^1(\mathbb{R})$, $g'(0)=0$, and there exist $C>0$ and $q\in(4,2\cdot2^*)$ with $2^*=2N/(N-2)$ such that
$|g'(t)|\leq C(1+|t|^{q-2})$.

\noindent($g_6$) $\lim_{|t|\rightarrow+\infty}\frac{g(t)}{t^3}=l\in (\mu,+\infty)$, where $\mu=\max\{\mu_1,\mu_2\}$ and $\mu_i$
is defined by
$$\mu_i=\inf\Bigl\{\int_{\Omega_i}|\nabla u|^2dx: u\in H^1_{0,rad}(\Omega_i), u>0\ \text{in}\ \Omega_i, \int_{\Omega_i}|u|^2dx=1\Bigr\}, \ i=1,2.$$

\begin{Remark}
An example of $V$ satisfying (V$_1$), (V$_3$) and (V$_4$) is
\begin{equation*}\aligned
V(x)=\left\{ \begin{array}{lll}
0,\ \quad\ & \text{if}\quad |x|\leq1,\\
(|x|-1)^2,\ \quad \ &\text{if}\quad
1<|x|\leq2,\\
1, \ \quad \ & \text{if}\quad |x|>2,
\end{array}\right.\endaligned
\end{equation*}
and in (V$_4$) we can choose $\Omega_1=\{x\in\mathbb{R}^N:|x|<\frac14\}$, $\Omega_2=\{x\in\mathbb{R}^N:\frac13<|x|<\frac12\}$.
\end{Remark}
\noindent\textbf{Theorem 1.3.} {\it Let (V$_1$), (V$_3$), (V$_4$), (g$_1$) and (g$_3$)-(g$_6$) hold. Then the problem (\ref{1.2}) has a least energy sign-changing solution, which changes sign exactly once.}

\begin{Remark}An interesting open problem left is whether the equation (\ref{1.2}) has multiple
sign-changing solutions if $\inf_{\mathbb{R}^N}V=0$.
\end{Remark}

 In the case $\inf_{\mathbb{R}^N}V>0$, Theorems 1.1 and 1.2 can be viewed
as the counterparts of the results of (\ref{1.2}) with 3-superlinear nonlinearity in \cite{DPW1,DPW2,YANG1} to
asymptotically 3-linear nonlinearity in some sense. Due to the competition of the quasilinear term $-u\Delta(u^2)$ with asymptotically 3-linear nonlinearity $g(u)$, we cannot use the methods of \cite{DPW1,DPW2,YANG1} trivially. We shall introduce a suitable restricted set, then we apply the method of sign-changing Nehari manifold dependent on the restricted set to show Theorems 1.1 and 1.2. In this process, on one hand, the  least energy of sign-changing  Nehari manifold has a good characterization which plays an important role in recovering the compactness of minimizing sequence. On the other hand,  the
involvement of the restricted set  in the method of sign-changing Nehari manifold enforces the implementation of some new estimates and verifications about the restricted set.

In the case $\inf_{\mathbb{R}^N}V=0$, the proof is quite different from the case $\inf_{\mathbb{R}^N}V>0$ and we will meet some new difficulties. Especially, we need to show that the sign-changing Nehari manifold is nonempty and the minimizing sequence on  the sign-changing Nehari manifold is a PS sequence. We shall use the solutions of two Dirichlet problems to construct a sign-changing function which belongs to the sign-changing Nehari manifold. Moreover, motivated by \cite{NW}, we shall apply Ekeland's variational principle and
implicit function theorem to obtain desired results. However, the quasilinear term $-u\Delta(u^2)$ leads that some new tricks and careful estimates of minimizing sequence are needed.

The paper is organized as
follows. In Section 2 we  give some preliminaries. In Sections 3, 4 and 5, we prove Theorems
1.1, 1.2 and 1.3 respectively.

\section{Preliminaries}
\setcounter{equation}{0}
In this paper we use the following
notations.  For $1\leq p\leq\infty$ and $\Omega\subset\mathbb{R}^N$, the
norms in $L^p(\mathbb{R}^N)$ and $L^p(\Omega)$ are denoted by $|\cdot|_{p}$ and $|\cdot|_{p,\Omega}$ respectively. $\int_{\mathbb{R}^N}
f(x)dx$ will be represented by $\int_{\mathbb{R}^N} f(x)$. For any $r>0$ and $x\in\mathbb{R}^N$,
$B_r(x)$ denotes the ball  centered at $x$ with the radius $r$.  For the set $\Lambda\subset\mathbb{R}^N$, $|\Lambda|$ denotes the Lebesgue measure of $\Lambda$. $H^1(\mathbb{R}^N)$ is the Hilbert space equipped with the standard norm
\begin{equation}\label{2.1}\|u\|^2_0=\int_{\mathbb{R}^N}(|\nabla u|^2+u^2).\end{equation}
The Hilbert space
$$E:= H^1_r(\mathbb{R}^N)=\bigl\{u\in H^1(\mathbb{R}^N):u(x)=u(|x|)\bigr\}$$ is the  space endowed with the following inner product and norm
$$(u,v)=\int_{\mathbb{R}^N} (\nabla u\nabla v+V(x)uv),\quad \|u\|^2=\int_{\mathbb{R}^N} (|\nabla u|^2+V(x)u^2).$$
Under (V$_1$) and (V$_2$), the norm $\|\cdot\|$ is an equivalent norm to the standard norm $\|\cdot\|_0$ in $E$,  and so the embedding $E\hookrightarrow L^p(\mathbb{R}^N)$ is compact for any $p\in (2,2^*)$.

 In view of  (V$_1$), (V$_3$), the H\"{o}lder and Sobolev inequalities we have
$$\aligned\int_{\mathbb{R}^N}(|\nabla u|^2+u^2)\leq& \int_{\mathbb{R}^N}|\nabla u|^2+|\mathcal{V}_a|^{\frac{2^*-2}{2^*}}\Bigl(\int_{\mathcal{V}_a}|u|^{2^*}\Bigr)^{\frac{2}{2^*}}
+a^{-1}\int_{\mathbb{R}^N\backslash{\mathcal{V}_a}}V(x)u^2\\ \leq& \max\{1+S^{-1}|\mathcal{V}_a|^{\frac{2^*-2}{2^*}}, a^{-1}\}\int_{\mathbb{R}^N}(|\nabla u|^2+V(x)u^2),\endaligned$$
where $S=\inf_{u\in H^1(\mathbb{R}^N)}|\nabla u|^2_2/{|u|^2_{2^*}}$.
On the other hand, since $V\in L^\infty(\mathbb{R}^N)$ in (V$_1$) we get $$\int_{\mathbb{R}^N}(|\nabla u|^2+V(x)u^2)\leq \max\{1,|V|_\infty\}\int_{\mathbb{R}^N}(|\nabla u|^2+u^2).$$  Hence, the norm $\|\cdot\|$ is an equivalent norm to the standard norm
 $\|\cdot\|_0$ in $E$.
 Then, under (V$_1$) and (V$_3$),  the embedding  $E\hookrightarrow L^p(\mathbb{R}^N)$ is compact for any $p\in (2,2^*)$.

Due to the term $-u\Delta(u^2)$, the functional of the problem (\ref{1.2})
$$J{(u)}=\frac{1}{2}\int_{\mathbb{R}^{N}} (1+2u^2)|\nabla u|^2+\frac 12\int_{\mathbb{R}^{N}}
V(x)u^2-\int_{\mathbb{R}^{N}} G(u),$$ is not well defined in $E$. As in
\cite{LWW}, we make use of a suitable change, namely, the change of
variables $v=f^{-1}(u)$, and then we can choose $E$
as the research space. The change $f$ is defined by $$\aligned
&f'(t)=\frac{1}{(1+2f^2(t))^{\frac12}} \ \text{on} \ [0,\infty),\\&
f(t)=-f(-t) \ \ \ \ \ \ \text{on} \ (-\infty,0].\endaligned$$ Below we state some properties of $f$.

\begin{lemma} \label{l2.1}
 \vskip 0.1 true cm \noindent(1) $f$ is uniquely defined, $C^\infty$
 and invertible;

\vskip 0.1 true cm \noindent (2) $|f'(t)|\leq1$ for all
$t\in\mathbb{R}$;

\vskip 0.1 true cm \noindent (3) $|f(t)|\leq|t|$ for all
$t\in\mathbb{R}$, and $|f(t)|\leq2^{\frac{1}{4}}|t|^{\frac12}$ for all $t\in\mathbb{R}$;

\vskip 0.1 true cm \noindent (4) $\frac{f(t)}{t}\rightarrow1$ as
$t\rightarrow0$;

\vskip 0.1 true cm \noindent (5)
$\frac{|f(t)|}{|t|^{\frac12}}\rightarrow2^{\frac{1}{4}}$ as
$|t|\rightarrow+\infty$;

\vskip 0.1 true cm \noindent (6) $\frac{f^2(t)}{2}\leq
tf'(t)f(t)\leq f^2(t)$ for all $t\in\mathbb{R}$;

\vskip 0.1 true cm \noindent (7)
the function $f(t)f'(t)t^{-1}$ is decreasing for all $t>0$,
and the function $f^3(t)f'(t)t^{-1}$ is  increasing for all $t>0$;

\vskip 0.2 true cm \noindent (8)
there exists a positive constant $C$ such that
$|f(t)|\geq
C|t|$ if $|t|\leq1$, and
$|f(t)|\geq
C|t|^\frac12$ if  $|t|\geq1$;
\vskip 0.2 true cm \noindent (9) $|f(t)f'(t)|\leq\frac{1}{\sqrt{2}}$ for all $t\in\mathbb{R}$;
\vskip 0.2 true cm \noindent (10) $f'(t)|t|^{\frac12}\rightarrow 2^{-\frac34}$ as $|t|\rightarrow+\infty$.
\end{lemma}
{\bf Proof}: The conclusions (1)-(9) were proved in \cite{CJ,DS1}. It suffices to show (10). Indeed, by the definition of $f'(t)$ and (5) we infer
$$\aligned\lim_{|t|\rightarrow+\infty}f'(t)|t|^{\frac12}
=\lim_{|t|\rightarrow+\infty}\frac{|t|^{\frac12}}{(1+2f^2(t))^{\frac12}}
=\lim_{|t|\rightarrow+\infty}\frac{|f(t)|}{(1+2f^2(t))^{\frac12}}\cdot\frac{|t|^{\frac12}}{|f(t)|}=2^{-\frac34}.\endaligned$$
Then (10) holds true. \ \ \  \ $\Box$

\begin{lemma}\label{l2.2}
(1)\ If ($g_1$) and ($g_2$) are satisfied, then for some $2<p<2^*$ and any $\epsilon>0$, there exists a constant
$C_{\epsilon}>0$ such that
\begin{equation}\label{2.2}|g(t)|\leq \epsilon|t|+C_\epsilon|t|^{2p-1} \ \ \text{and}\ |G(t)|\leq\frac{\epsilon}{2}t^2+\frac{C_\epsilon}{p}|t|^{2p},\ \forall t\in\mathbb{R},\end{equation}
and
\begin{equation}\label{2.3}|g(f(t))f(t)|\leq \epsilon t^2+C_\epsilon|t|^{p} \ \ \text{and}\ |G(f(t))|\leq\frac{\epsilon}{2}t^2+\frac{C_\epsilon}{p}|t|^{p},\ \forall t\in\mathbb{R}.\end{equation}
(2) If ($g_2$) and ($g_3$) are satisfied, then $\frac{g(t)}{t^3}<1$ for all $t\in\mathbb{R}\backslash\{0\}$.\\
(3) If ($g_1$) and ($g_3$) are satisfied, then
$0\leq4G(t)\leq g(t)t$ for any $t\in\mathbb{R}$.
\end{lemma}

By the definition and properties of the change $f$, from ${J}$ we obtain
the functional
\begin{equation*}I(v)=\frac12\int_{\mathbb{R}^N}|\nabla v|^2+\frac12\int_{\mathbb{R}^N}
V(x)f^2(v)-\int_{\mathbb{R}^N} G(f(v)),\end{equation*} which is well defined in
$E$ and of $C^1$ under our hypotheses. Moreover, the
critical points of $I$ are the weak solutions of the
problem\begin{equation}\label{2.2.0}
-\Delta v+V( x)f(v)f'(v)=g(f(v))f'(v),\ \ x\in\mathbb{R}^N.
\end{equation}
By the one-to-one correspondence $f$, we just need to study the equation (\ref{2.2.0}).

\begin{lemma}\label{l2.3} Let (V$_1$), and (V$_2$) or (V$_3$) hold.
 Then there exists $C>0$ such that
 \begin{equation}\label{2.3.4}\|v\|^2\leq C\Bigl[|\nabla v|^{2}_2+|\nabla v|^{2^*}_2+\int_{\mathbb{R}^N}V(x)f^2(v)\Bigr],\ \ \ \forall v\in E.\end{equation}
\end{lemma}
{\bf Proof}: Observe that
\begin{equation*}\int_{\{|v|>1\}}V(x)v^2\leq |V|_\infty\int_{\{|v|>1\}}|v|^{2^*}
\leq C|V|_\infty|\nabla v|^{2^*}_2.\end{equation*}
Moreover, by Lemma 2.1 (8) we have
\begin{equation*}\int_{\{|v|\leq1\}}V(x)v^2\leq\frac1 {C^2}\int_{\{|v|\leq1\}}V(x)f^2(v)\leq\frac1 {C^2}\int_{\mathbb{R}^N}V(x)f^2(v).\end{equation*}
Then (\ref{2.3.4}) yields.\ \ \ $\Box$

As the proof of (3.2) in \cite{Fang}, we know
\begin{lemma}\label{l2.4}Under the assumptions of Lemma \ref{l2.3}, there exist $C_1,\rho>0$ such that
\begin{equation}\label{2.4}
|\nabla u|^2_2+\int_{\mathbb{R}^N}V(x)f^2(u)\geq C_1\|u\|^2, \ \text{whenever}\ \|u\|\leq\rho.
\end{equation}\end{lemma}

Next we define the Nehari manifold associated to $I$ by
$$\mathcal{N}=\Bigl\{u\in E\backslash\{0\}:\langle I'(u),u\rangle=0\Bigr\},$$
where $$\langle I'(u),u\rangle=|\nabla u|^2_2+\int_{\mathbb{R}^N}V(x)f(u)f'(u)u-\int_{\mathbb{R}^N} g(f(u))f'(u)u,$$
and the least energy on $\mathcal{N}$ is defined as \begin{equation}\label{2.5}d=\inf_{u\in \mathcal{N}}I(u).\end{equation}
We also define the sign-changing Nehari manifold associated to $I$ by
$$\mathcal{M}=\{u\in E: u^{\pm}\neq0\ \text{and}\ \langle I'(u),u^{\pm}\rangle=0\},$$
where
$$\langle I'(u),u^{\pm}\rangle=|\nabla u^{\pm}|^2_2+\int_{\mathbb{R}^N}V(x)f(u^{\pm})f'(u^{\pm})u^{\pm}-\int_{\mathbb{R}^N} g(f(u^{\pm}))f'(u^{\pm})u^{\pm},$$
and the least energy on $\mathcal{M}$ is defined as
\begin{equation*}c=\inf_{u\in \mathcal{M}}I(u).\end{equation*}

\section{Proof of Theorem 1.1}

In this section, we are devoted to showing Theorem 1.1 and always assume (V$_1$), (V$_2$) and (g$_1$)-(g$_4$) are satisfied. To adapt the method of Nehari manifold, we introduce a set
$$\Theta=\{u\in E: |\nabla u^{\pm}|^2_2<|u^{\pm}|^2_2\}.$$
It  is easy to see that $\Theta\neq\emptyset$. Moreover, by the properties of the change $f$, we have $\mathcal{M}\subset\Theta$.
\begin{lemma}\label{l3.1}For any $u\in \Theta$, there exists a unique pair $(s_u,t_u)$ of positive numbers
such that $s_uu^++t_uu^-\in \mathcal{M}$ and $I(s_uu^++t_uu^-)=\max_{s,t>0}I(su^++tu^-)$. Moreover,
\begin{equation}\label{3.1.0}c=\inf_{\mathcal{M}}I=\inf_{u\in \Theta}\max_{s,t>0}I(su^++tu^-).\end{equation}\end{lemma}
{\bf Proof}: Set $\psi_1(s,t)=\langle I'(su^++tu^-),su^+\rangle$, $\psi_2(s,t)=\langle I'(su^++tu^-),tu^-\rangle$ for $s,t>0$, namely
$$\psi_1(s,t)=s^2|\nabla u^+|^2_2+\int_{\mathbb{R}^N}V(x)f(su^+)f'(su^+)su^+-\int_{\mathbb{R}^N}g(f(su^+))f'(su^+)su^+,$$
and
$$\psi_2(s,t)=t^2|\nabla u^-|^2_2+\int_{\mathbb{R}^N}V(x)f(tu^-)f'(tu^-)tu^--\int_{\mathbb{R}^N}g(f(tu^-))f'(tu^-)tu^-.$$
By (\ref{2.3}) and Lemma \ref{l2.1} (6) we know
$$\psi_1(s,t)\geq s^2|\nabla u^+|^2_2+\int_{\mathbb{R}^N}V(x)f(su^+)f'(su^+)su^+-\epsilon s^2|u^+|^2_2-C_\epsilon s^p|u^+|^p_p,$$
for some $p\in(2,2^*)$ and any $\epsilon>0$.
Using (\ref{2.4}) and Lemma \ref{l2.1} (6) we get
$$s^2|\nabla u^+|^2_2+\int_{\mathbb{R}^N}V(x)f(su^+)f'(su^+)su^+\geq C_1s^2\|u^+\|^2,$$ if small  $s$  satisfies $s^2\|u^+\|\leq\rho$. Then for small $s$ satisfying $s^2\|u^+\|\leq\rho$, we have
$$\psi_1(s,t)\geq C_1s^2\| u^+\|^2-C\epsilon s^2\|u^+\|^2-C_\epsilon s^p\|u^+\|^p.$$
Choosing sufficiently small $\epsilon>0$, there exists small $s_m>0$ such that
\begin{equation}\label{2.2.2}\psi_1(s_m,t)>0, \ \text{for all}\ t>0.\end{equation} Similarly, there exists small $t_m>0$ such that $\psi_2(s,t_m)>0$ for all $s>0$.

On the other hand, note that
$$\frac{\psi_1(s,t)}{s^2}=|\nabla u^+|^2_2+\int_{\mathbb{R}^N}V(x)\frac{f(su^+)f'(su^+)su^+}{s^2}
-\int_{\mathbb{R}^N}\frac{g(f(su^+))f'(su^+)su^+}{s^2}.$$
By Lemma 2.1 (5) and (6) we have
$$\lim_{s\rightarrow+\infty}\int_{\mathbb{R}^N}V(x)\frac{f(su^+)f'(su^+)su^+}{s^2}=0.$$
In view of (g$_2$), Lemma \ref{l2.1} (5) and (10) we deduce
\begin{equation}\label{3.4}\aligned&\lim_{s\rightarrow+\infty}\int_{\mathbb{R}^N}\frac{g(f(su^+))f'(su^+)su^+}{s^2}
\\=&\lim_{s\rightarrow+\infty}\int_{\mathbb{R}^N}
\frac{g(f(su^+))}{f^3(su^+)}\cdot \frac{f^3(su^+)}{s^{\frac32}(u^+)^{\frac32}}\cdot f'(su^+)s^{\frac12}(u^+)^{\frac12}\cdot(u^+)^{2}=|u^+|^2_2.\endaligned\end{equation}
Since $u\in \Theta$, we have $|\nabla u^+|^2_2<|u^+|^2_2$. So $$\lim_{s\rightarrow+\infty}\frac{\psi_1(s,t)}{s^2}=|\nabla u^+|^2_2-|u^+|^2_2<0.$$ Then there exists $s_M>0$ such that $\psi_1(s_M,t)<0$ for all $t>0$. Similarly, there exists $t_M>0$ such that $\psi_2(s,t_M)<0$ for all $s>0$.

Hence, by taking $r=\min\{t_m,s_m\}$ and $R=\max\{t_M,s_M\}$, we infer
$$\psi_1(r,t)>0\ \text{and}\ \psi_1(R,t)<0\ \ \text{for all}\ t\in[r,R],$$
and
$$\psi_2(s,r)>0\ \text{and}\ \psi_2(s,R)<0\ \ \text{for all}\ s\in[r,R].$$
Applying \cite[Lemma 2.4]{LiLuo}, there exists $(s_u,t_u)\in(r,R)\times(r,R)$ such that $\psi_1(s_u,t_u)=\psi_2(s_u,t_u)=0$. That is, $s_uu^++t_uu^-\in \mathcal{M}$.

Below we claim that such $(s_u,t_u)$ is unique. We firstly show that if $u\in \mathcal{M}$, then $(s_u,t_u)=(1,1)$. Indeed, since $s_uu^++t_uu^-\in \mathcal{M}$, we have
\begin{equation}\label{2.1.6}|\nabla u^+|^2_2=-\int_{\mathbb{R}^N}V(x)\frac{f(s_uu^+)f'(s_uu^+)s_uu^+}{s^2_u}+\int_{\mathbb{R}^N}\frac{g(f(s_uu^+))f'(s_uu^+)s_uu^+}{s^2_u},\end{equation}
and \begin{equation}\label{2.1.7}|\nabla u^-|^2_2=-\int_{\mathbb{R}^N}V(x)\frac{f(t_uu^-)f'(t_uu^-)t_uu^-}{t^2_u}+\int_{\mathbb{R}^N}\frac{g(f(t_uu^-))f'(t_uu^-)t_uu^-}{t^2_u}.\end{equation}
Note that $u\in \mathcal{M}$, then
\begin{equation}\label{2.1.8}|\nabla u^+|^2_2+\int_{\mathbb{R}^N}V(x)f(u^+)f'(u^+)u^+=\int_{\mathbb{R}^N}g(f(u^+))f'(u^+)u^+,\end{equation}
and
\begin{equation}\label{2.1.9}|\nabla u^-|^2_2+\int_{\mathbb{R}^N}V(x)f(u^-)f'(u^-)u^-=
\int_{\mathbb{R}^N}g(f(u^-))f'(u^-)u^-.\end{equation}
By Lemma 2.1 (7) and ($g_3$), we know \begin{equation}\label{3.7.0}-\int_{\mathbb{R}^N}V(x)\frac{f(su^+)f'(su^+)su^+}{s^2}
+\int_{\mathbb{R}^N}\frac{g(f(su^+))f'(su^+)su^+}{s^2}\ \text{is increasing in} \ s.\end{equation}
Then from (\ref{2.1.6}) and (\ref{2.1.8}) we get $s_u=1$. Similarly, from (\ref{2.1.7}) and (\ref{2.1.9}) we have $t_u=1$.

Next we show that if $u\in \Theta\backslash \mathcal{M}$, then $(s_u,t_u)$ is unique.
Suppose that there exists another pair $(s'_u,t'_u)$ of positive numbers such that $s'_uu^++t'_uu^-\in \mathcal{M}$. Then
$$\frac{s'_u}{s_u}(s_uu^+)+\frac{t'_u}{t_u}(t_uu^-)\in \mathcal{M}.$$
Therefore, $s'_u=s_u$ and $t'_u=t_u$. Namely, $(s_u,t_u)$ is unique.

Define $\Psi:\mathbb{R}^+\times\mathbb{R}^+\rightarrow\mathbb{R}$ by $\Psi(s,t)=I(su^++tu^-)$. Then we know $(s_u, t_u)$ is the unique critical point of $\Psi$ in $\mathbb{R}^+\times\mathbb{R}^+$. It is easy to see that $I(su^++tu^-)\rightarrow-\infty$ as $|(s,t)|\rightarrow+\infty$. If we assume that $(0,t_0)$ is a  maximum point of $\Psi$ with $t_0\geq0$, as (\ref{2.2.2}) we infer
$$\frac{\partial \Psi(s,t_0)}{\partial s}=s^{-1}\psi_1(s,t_0)>0,$$
for small enough $s$. Therefore, $\Psi(s,t_0)$ is an increasing function with respect to $s$ if $s$ is small enough. Then the pair $(0,t_0)$ is not a maximum point of $\Psi$  in $\mathbb{R}^+\times\mathbb{R}^+$.
This ends the proof.\ \ \ \ \ $\Box$



\begin{lemma}\label{l3.4}(1) There exists $r>0$ such that  $\|u^{\pm}\|\geq r$ for any $u\in \mathcal{M}$, and $\|u\|\geq r$ for any $u\in \mathcal{N}$;\\ (2) There exists $\bar{r}>0$ such that\ $c=\inf_{\mathcal{M}}I\geq \bar{r}>0$.\end{lemma}
{\bf Proof}: (1)
Set $u\in \mathcal{M}$. For any $\epsilon>0$ and some $p\in (2,2^*)$, by (\ref{2.2}), Lemma \ref{l2.1} (3) and (6), there holds
\begin{equation}\label{2.2.1}\aligned
&|\nabla u^{+}|^2_2+\int_{\mathbb{R}^N}V(x)f(u^{+})f'(u^{+})u^{+}=\int_{\mathbb{R}^N}g(f(u^{+}))f'(u^+)u^+\\
&\leq \epsilon|f(u^+)|^2_2+C_\epsilon|f(u^+)|^{2p}_{2p}\leq C\epsilon\|u^+\|^2+C_\epsilon\|u^{+}\|^{p}.
\endaligned\end{equation}
If $\|u^+\|\geq\rho$ with $\rho$ given in (\ref{2.4}), we are done. Otherwise, if $\|u^+\|\leq\rho$, by (\ref{2.4}) and Lemma \ref{l2.1} (6) we then infer
$$C_1\|u^+\|^2\leq C\epsilon\|u^+\|^2+C_\epsilon\|u^+\|^p.$$
Then choosing $\epsilon$ small and $r<\rho$ we get
$$\|u^{+}\|\geq r,\ \forall u\in \mathcal{M}.$$
Similarly, we have $\|u^{-}\|\geq r,$ and $\|u\|\geq r$ for all $u\in \mathcal{N}$.

(2) By (\ref{2.3}) we know
$$I(u^+)\geq\frac12\Bigl(|\nabla u^+|^2_2+\int_{\mathbb{R}^N}V(x)f^2(u^+)\Bigr)
-\epsilon|u^+|^2_2-C_\epsilon|u^+|^p_p,$$
for some $p\in (2,2^*)$ and any $\epsilon>0$. Using (\ref{2.4}) we get
$$I(u^+)\geq\frac{C_1}{2}\|u^+\|^2-\epsilon|u^+|^2_2-C_\epsilon|u^+|^p_p,$$
if $\|u^+\|\leq\rho$. Choosing small $\epsilon>0$ and $\bar{\rho}<\rho$, we infer
\begin{equation}\label{3.8}\inf_{S^+_{\bar{\rho}}}I\geq\bar{r}>0, \ \text{where}\ S^+_{\bar{\rho}}:=\{u\in E:\|u^+\|=\bar{\rho}\}.\end{equation}
Since $u\in \mathcal{M}\subset\Theta$, (\ref{3.8}) and (\ref{3.1.0})
imply that $c=\inf_{\mathcal{M}}I\geq\bar{r}>0$. \ \ \ \ \ $\Box$

\begin{lemma}\label{l3.3} There is $\varrho>0$ such that
\begin{equation}\label{3.1.1}\int_{\mathbb{R}^N}(|\nabla u^{\pm}|^2+V(x)f^2(u^{\pm}))\geq\varrho,\ \ \forall u\in\mathcal{M},\end{equation}
and
\begin{equation}\label{3.3}|f(u^{\pm})|^{2p}_{2p}\geq\varrho,\ \ \forall u\in\mathcal{M},\end{equation}
where $p\in (2,2^*)$.\end{lemma}
{\bf Proof}: Since $u\in \mathcal{M}$, by the first inequality of (\ref{2.2.1}) and (V$_2$) we get
\begin{equation}\label{3.5}|\nabla u^{+}|^2_2+\int_{\mathbb{R}^N}V(x)f^2(u^+)\leq C|f(u^+)|^{2p}_{2p}.\end{equation}
Note that $2p\in (4,22^*)$, we may assume that $2p=2\tau+22^*(1-\tau)$, $\tau\in (0,1)$. Then
$$\aligned |f(u^+)|^{2p}_{2p}\leq |f(u^+)|^{2\tau}_2|f(u^+)|^{22^*(1-\tau)}_{2^*}&\leq |f(u^+)|^{2\tau}_2|u^+|^{2^*(1-\tau)}_{2^*}\leq |f(u^+)|^{2\tau}_2|\nabla u^+|^{2^*(1-\tau)}_{2}\\&\leq C\Bigl(\int_{\mathbb{R}^N}V(x)f^2(u^+)+|\nabla u^+|^2_2\Bigr)^{\tau+\frac{2^*}{2}(1-\tau)}.\endaligned$$
Since $\tau+\frac{2^*}{2}(1-\tau)>1$, (\ref{3.1.1}) holds true.
In view of (\ref{3.1.1}) and (\ref{3.5}), (\ref{3.3}) follows. \ \ \ \ $\Box$

{\bf Proof of Theorem 1.1.}
Assume that $\{u_n\}\subset \mathcal{M}$ is a minimizing sequence, i.e. $I(u_n)\rightarrow c$.

{\bf Claim 1}: $\{u_n\}$ is  bounded in $E$.

In fact, arguing by contradiction we assume that $\|u_n\|\rightarrow\infty$. Set $v_n=\frac{u_n}{\|u_n\|}$. Then up to a subsequence, we may suppose $v_n\rightharpoonup v$ in $E$, $v_n\rightarrow v$ in $L^r(\mathbb{R}^N)$ with $2<r<2^*$ and $v_n(x)\rightarrow v(x)$
a.e. in $\mathbb{R}^N$. If $|v_n|_r\rightarrow0$, by (\ref{2.3}) we obtain $\int_{\mathbb{R}^N}G(f(sv_n))\rightarrow0$ for all $s>0$. From Lemma \ref{l2.1} (8) and Lemma \ref{l3.1} it follows that
\begin{equation}\label{3.12.0}\aligned
c+o_n(1)&\geq I(u_n)\geq I(sv_n)
\\&=\frac{s^2}{2}|\nabla v_n|^2_2+\frac{1}{2}\int_{\mathbb{R}^N}V(x)f^2(sv_n)-\int_{\mathbb{R}^N}G(f(sv_n))\\
&\geq\frac{s^2}{2}|\nabla v_n|^2_2+\frac{C^2s^2}{2}\int_{\{|sv_n|\leq1\}}V(x)v^2_n
-\int_{\mathbb{R}^N}G(f(sv_n))\\
&=\frac{s^2}{2}|\nabla v_n|^2_2+\frac{C^2s^2}{2}\int_{\mathbb{R}^N}V(x)v^2_n-\frac{C^2}{2}\int_{\{|sv_n|\geq1\}}V(x)(sv_n)^2
-\int_{\mathbb{R}^N}G(f(sv_n))\\
&\geq\frac{s^2}{2}\min\{1,C^2\}-C_1\int_{\mathbb{R}^N}|sv_n|^{p}-\int_{\mathbb{R}^N}G(f(sv_n))
\rightarrow\frac{s^2}{2}\min\{1,C^2\},
\endaligned\end{equation}
where $p\in(2,2^*)$. Taking sufficiently large $s$, there is a contradiction. Hence, $|v_n|_r\not\rightarrow0$. From \cite[Lemma 1.21]{WM1} it follows that  there exist $\delta>0$ and $y_n\in\mathbb{R}^N$ such that
$\int_{B_{1}(y_n)}v^2_n(x)dx>\delta.$
Set $w_n(x)=v_n(x+y_n)$. Then there exists $w_0\neq0$ such that $w_n\rightharpoonup w_0$ in $E$. Let $\Lambda:=\{x\in\mathbb{R}^N: w_0(x)\neq0\}$. Therefore $|u_n(x+y_n)|\rightarrow\infty$ for $x\in \Lambda$. By ($g_4$) and Lemma \ref{l2.2} (3) we deduce
\begin{equation}\label{3.10.0}
\aligned
&c+o_n(1)=I(u_n)-\frac12\langle I'(u_n),u_n\rangle\\
&=\frac12\int_{\mathbb{R}^N}V(x)\bigl[f^2(u_n)-f'(u_n)f(u_n)u_n\bigr]+
\int_{\mathbb{R}^N}\bigl[\frac12g(f(u_n))f'(u_n)u_n-G(f(u_n))\bigr]\\
&\geq\int_{\mathbb{R}^N}\bigl[\frac14g(f(u_n))f(u_n)-G(f(u_n))\bigr]\\
&\geq\int_{\Lambda}\bigl[\frac14g(f(u_n(x+y_n)))f(u_n(x+y_n))-G(f(u_n(x+y_n)))\bigr]\rightarrow+\infty.
\endaligned\end{equation}
This is a contradiction.
Therefore, $\{u_n\}$ is bounded in $E$ and assume that $u_n\rightharpoonup u$ in $E$. Then it is easy to see that
$u^+_n\rightharpoonup u^+$  and $u^-_n\rightharpoonup u^-$ in $E$. Up to a subsequence, we suppose that $u^\pm_n\rightarrow u^\pm$ in $L^p(\mathbb{R}^N)$, $2<p<2^*$. Then
$f(u^\pm_n)\rightarrow f(u^\pm)$ in $L^r(\mathbb{R}^N)$, $2<r<22^*$.
By (\ref{3.3}), we know $u^\pm\neq0$. It is easy to see that
$$\int_{\mathbb{R}^N}g(f(u^\pm_n))f'(u^\pm_n)u^\pm_n\rightarrow\int_{\mathbb{R}^N}g(f(u^\pm))f'(u^\pm)u^\pm.$$
From Fatou lemma it follows that
\begin{equation}\label{zh2}\aligned\langle &I'(u^\pm),u^\pm\rangle=|\nabla u^\pm|^2_2+\int_{\mathbb{R}^N}V(x)f(u^\pm)f'(u^\pm)u^\pm-
\int_{\mathbb{R}^N}g(f(u^\pm))f'(u^\pm)u^\pm\\&\leq\liminf_{n\rightarrow\infty}\Bigl[|\nabla u^\pm_n|^2_2+\int_{\mathbb{R}^N}V(x)f(u^\pm_n)f'(u^\pm_n)u^\pm_n\Bigr]-
\lim_{n\rightarrow\infty}\int_{\mathbb{R}^N}g(f(u^\pm_n))f'(u^\pm_n)u^\pm_n\\
&=\liminf_{n\rightarrow\infty}\langle I'(u^\pm_n),u^\pm_n\rangle=0.\endaligned\end{equation}  In view of $u^\pm\neq0$, we have $\int_{\mathbb{R}^N}f(u^{\pm})f'(u^{\pm})u^{\pm}>0$.  Then from (\ref{zh2}), Lemma \ref{l2.1} (3) and (9) we obtain
\begin{equation}\label{3.14.0}\aligned|\nabla u^{\pm}|^2_2&<\int_{\mathbb{R}^N}g(f(u^{\pm}))f'(u^{\pm})u^{\pm}
<\int_{\mathbb{R}^N}f^3(u^{\pm})f'(u^{\pm})u^{\pm}\leq|u^{\pm}|^2_2.
\endaligned\end{equation} So $u\in \Theta$. By Lemma \ref{l3.1} there exist $\alpha,\beta>0$ such that $\alpha u^++\beta u^-\in\mathcal{M}$.

{\bf Claim 2}: $c$ is attained and $I(\alpha u^++\beta u^-)=c$.

Indeed, one easily has that $I$ is weak lower semicontinuous. Then
\begin{equation}\label{3.18.0}c\leq I(\alpha u^++\beta u^-)\leq\liminf_{n\rightarrow\infty}I(\alpha u^+_n+\beta u^-_n)\leq\liminf_{n\rightarrow\infty}I(u_n)=c.\end{equation}
Therefore, $I(\alpha u^++\beta u^-)=c$.

\textbf{Claim 3}: $I'(w)=0$, where $w=\alpha u^++\beta u^-\in \mathcal{M}$ obtained above.

Indeed, on the contrary, we assume $I'(w)\neq0$. Then there exist $\delta,\nu>0$ such that
$$\|I'(v)\|\geq\nu \ \ \text{for all} \ v\in E\ \text{satisfying}\ \|v-w\|\leq 2\delta.$$
Since $w^+\neq0$ and $w^-\neq0$, we take $L=\min\bigl\{{(\inf_{\mathbb{R}^N}V)}^\frac12|w^+|_2,{(\inf_{\mathbb{R}^N}V)}^\frac12|w^-|_2\bigr\}$, and without loss of generality,
we may assume $6\delta<L$.

Let $D=[\frac12,\frac32]\times[\frac12,\frac32]$. Since $w\in \mathcal{M}\subset \Theta$, from Lemma \ref{l3.1} we know
$$I(tw^++sw^-)<I(w)=c,$$
holds for $(t,s)\in D$ with $t\neq1$ or $s\neq1$. Then
\begin{equation}\label{2.4.1} \bar{c}=\max_{\partial D}I(tw^++sw^-)<c.\end{equation}
Applying the quantitative deformation lemma \cite[Lemma 2.3]{WM1} with $\epsilon=\min\{\frac{(c-\bar{c})}{4},\frac{\nu\delta}{8}\}$ and $S=B(w,\delta)$, there exists $\eta\in C([0,1]\times E, E)$ such that

\noindent(a) $\eta(1,v)=v$ if $v\not\in I^{-1}([c-2\epsilon,c+2\epsilon])\cap S_{2\delta}$,\\
\noindent(b) $\eta(1,I^{c+\epsilon}\cap S)\subset I^{c-\epsilon}$, where $I^{c\pm\epsilon}:=\{v\in E: I(v)\leq c\pm\epsilon\}$,\\
\noindent(c) $I(\eta(1,v))\leq I(v)$ for all $v\in E$,
\\
\noindent(d) $\|\eta(t,v)-v\|\leq\delta$ for all $t\in[0,1]$ and $v\in E$.

It is clear that
\begin{equation}\label{2.4.2}\max_{(s,t)\in D}I(\eta(1,tw^++sw^-))<c.\end{equation}
Now we claim that there exists $(t_0,s_0)\in D$ such that
\begin{equation}\label{2.4.3}\eta(1,t_0w^++s_0w^-)\in \mathcal{M}.\end{equation}
In fact, define $\varphi(t,s)=\eta(1,tw^++sw^-)$ and
$$\Psi_0(t,s)=\bigl(\langle I'(tw^++sw^-),tw^+\rangle, \langle I'(tw^++sw^-),sw^-\rangle\bigr).$$$$\Psi_1(t,s)=\bigl(\langle I'(\varphi(t,s)),\varphi^+(t,s)\rangle, \langle I'(\varphi(t,s)),\varphi^-(t,s)\rangle\bigr).$$
By Lemma \ref{l3.1} and the degree theory, we get $deg(\Psi_0, D,0)=1$. From (\ref{2.4.1}) and (a) it follows that $\varphi(t,s)=tw^++sw^-$ on $\partial D$. Thus,
using the degree theory, we have $deg(\Psi_1, D,0)=deg(\Psi_0, D,0)=1$. Consequently, $\Psi_1(t_0,s_0)=(0,0)$ for some $(t_0,s_0)\in D$. If $[\varphi(t_0,s_0)]^+=0$, by the Young inequality we have
$$\aligned&\|\varphi(t_0,s_0)-(t_0w^++s_0w^-)\|\\ \geq& \Bigl(\int_{\mathbb{R}^N}V(x)|t_0w^++s_0w^--[\varphi(t_0,s_0)]^-|^2\Bigr)^{\frac12}\geq (\inf_{\mathbb{R}^N}V)^{\frac12}t_0|w^+|_2>2\delta,\endaligned$$
which contradicts with (d). So
\begin{equation}\label{3.17.0}[\varphi(t_0,s_0)]^+\neq0.\end{equation} Similarly, we can show $[\varphi(t_0,s_0)]^-\neq0$. Therefore, $\eta(1,t_0w^++s_0w^- )=\varphi(t_0,s_0)\in \mathcal{M}$. This is a contradiction with (\ref{2.4.2}).

Now following the arguments
in \cite{Chentang} (or see \cite[p. 1051]{CCN}) we show that the solution $w$ changes sign just once in $\mathbb{R}^N$. Assume that $w=w_1+w_2+w_3$ satisfies
$$w_i\neq0,\ w_1\geq0, \ w_2\leq0,\ \text{and}\ supp(w_i)\cap supp(w_j)=\emptyset, \ \text{for}\ i\neq j,\ \text{and}\ i,j=1,2,3.$$
Using $I'(w)=0$ we know $\langle I'(w_i),w_i\rangle=0$ for $i=1,2,3$. Note that
$I(w_3)=I(w_3)-\frac12\langle I'(w_3),w_3\rangle>0$, we infer
$$I(w_1+w_2)<I(w_1+w_2+w_3)=I(w)=c,$$
a contradiction since $w_1+w_2\in \mathcal{M}$. This ends the proof.\ \ \ \ $\Box$
\section{Proof of Theorem 1.2}
In this section, we will prove Theorem 1.2 and always assume (V$_1$), (V$_2$) and (g$_1$)-(g$_4$) are satisfied. To get the solutions of (\ref{2.2.0}) possessing multiple nodal domains, we firstly look for positive solutions of (\ref{2.2.0}) on rotationally symmetric domains.
\subsection{The existence of positive solutions}
For any  rotationally symmetric domain $\Omega\subset\mathbb{R}^N$, we denote $H^1_{0,rad}(\Omega)$ as the Hilbert space of radial functions in $H^1_0(\Omega)$ with the following inner product  and norm  $$(u,v)_\Omega=\int_{\Omega}(\nabla u \nabla v+V(x)uv),\quad \|u\|_\Omega=\int_{\Omega}(|\nabla u|^2+V(x)u^2).$$
For $0\leq\rho<\sigma\leq+\infty$, define
$$\Omega(\rho,\sigma)=int\{x\in\mathbb{R}^N:\rho\leq|x|<\sigma\}.$$
In what follows, we consider the problem
\begin{equation}\label{5.1}\aligned
\left\{ \begin{array}{lll}
-\Delta u+V(x)f(u)f'(u)=g(f(u))f'(u),\ & \text{in}\quad \Omega(\rho,\sigma),\\
u\in E_{\rho,\sigma},
\end{array}\right.\endaligned
\end{equation}
where $E_{\rho,\sigma}:=H^1_{0,rad}(\Omega(\rho,\sigma))$.
 The functional $I_{\rho,\sigma}: E_{\rho,\sigma}\rightarrow \mathbb{R}$ is defined by
$$I_{\rho,\sigma}(u)=\frac12\int_{\Omega(\rho,\sigma)}(|\nabla u|^2+V(x)f^2(u))-\int_{\Omega(\rho,\sigma)}G(f(u)).$$
The corresponding Nehari manifold is
$$\mathcal{N}_{\rho,\sigma}:=\Bigl\{u\in E_{\rho,\sigma}\backslash\{0\}: \langle I'_{\rho,\sigma}(u),u\rangle=0\Bigr\},$$
and the least energy is
$$d_{\rho,\sigma}:=\inf_{\mathcal{N}_{\rho,\sigma}}I_{\rho,\sigma}.$$
Moreover, in order to use the method of Nehari manifold, we set
$$\Theta_{\rho,\sigma}:=\bigl\{u\in E_{\rho,\sigma}:|\nabla u|^2_{2,\Omega(\rho,\sigma)}<|u|^2_{2,\Omega(\rho,\sigma)}\bigr\}.$$
We understand that
$E_{\rho,\sigma}\subset E$ and $\mathcal{N}_{\rho,\sigma}\subset\mathcal{N}$ by defining $u(x)=0$ for $x\not\in \Omega(\rho,\sigma)$. As in Section 3, we have the following results.
\begin{lemma}\label{l5.1.1} (i) For any $u\in \Theta_{\rho,\sigma}$, there exists a unique $t_u>0$
such that $t_uu\in \mathcal{N}_{\rho,\sigma}$ and $I_{\rho,\sigma}(t_uu)=\max_{t>0}I_{\rho,\sigma}(tu)$. Moreover,
$d_{\rho,\sigma}=\inf_{u\in \Theta_{\rho,\sigma}}\max_{t>0}I_{\rho,\sigma}(tu);$

\noindent(ii) $\|u\|_{\rho,\sigma}\geq r>0$ for any $u\in \mathcal{N}_{\rho,\sigma}$ and $d_{\rho,\sigma}\geq \bar{r}>0$;

\noindent(iii) There is $\varrho>0$ such that
\begin{equation}\label{5.1.3}\int_{\Omega(\rho,\sigma)}|f(u)|^{2p}\geq\varrho,\ \ \forall u\in\mathcal{N}_{\rho,\sigma},\end{equation}
where $p\in(2,2^*)$.\end{lemma}

Below we are devoted to showing that $d_{\rho,\sigma}$ is achieved and the corresponding minimizer is a solution of the equation (\ref{5.1}).

\begin{lemma}\label{l5.1.2} $d_{\rho,\sigma}=\inf_{ \mathcal{N}_{\rho,\sigma}}I_{\rho,\sigma}$ is achieved by a positive function $w$ which is a solution of the equation (\ref{5.1}).\end{lemma}
{\bf Proof}: Assume that $\{u_n\}\subset \mathcal{N}_{\rho,\sigma}$ satisfies $I_{\rho,\sigma}(u_n)\rightarrow d_{\rho,\sigma}$. We firstly show that $\{u_n\}$ is  bounded in $E_{\rho,\sigma}$.
In fact, arguing by contradiction we assume that $\|u_n\|_{\Omega(\rho,\sigma)}\rightarrow+\infty$. Set $v_n=\frac{u_n}{\|u_n\|_{\Omega(\rho,\sigma)}}$. Up to a subsequence, we suppose $v_n\rightharpoonup v$ in $E_{\rho,\sigma}$, $v_n\rightarrow v$ in $L^r(\Omega(\rho,\sigma))$ with $2<r<2^*$ and $v_n(x)\rightarrow v(x)$
a.e. in $\Omega(\rho,\sigma)$. If $v=0$, then by (\ref{2.3}) we obtain $\int_{\Omega(\rho,\sigma)}G(f(sv_n))\rightarrow0$ for all $s>0$. Using Lemma \ref{l5.1.1}, similar to (\ref{3.12.0}) we get
$$\aligned
d_{\rho,\sigma}+o_n(1)&\geq I_{\rho,\sigma}(u_n)\geq I_{\rho,\sigma}(sv_n)
\rightarrow\frac{s^2}{2}\min\{1,C^2\}.
\endaligned$$
 Taking sufficiently large $s$, there is a contradiction. Hence, $v\neq0$ in $\Omega(\rho,\sigma)$. As (\ref{3.10.0}) we deduce
\begin{equation*}\aligned d_{\rho,\sigma}+o_n(1)&=I_{\rho,\sigma}(u_n)
\geq\int_{\Omega(\rho,\sigma)}\bigl[\frac14g(f(u_n))f(u_n)-G(f(u_n))\bigr]
\rightarrow+\infty.
\endaligned\end{equation*}
This is a contradiction.
Therefore, $\{u_n\}$ is bounded in $E_{\rho,\sigma}$ and assume that $u_n\rightharpoonup u$ in $E_{\rho,\sigma}$ and $u_n\rightarrow u$ in $L^r(\Omega(\rho,\sigma))$ with $2<r<2^*$. Then
$f(u_n)\rightarrow f(u)$ in $L^r(\Omega(\rho,\sigma))$, $2<r<22^*$.
By Lemma \ref{l5.1.1} (iii), we know $u\neq0$ in $\Omega(\rho,\sigma)$. Moreover, as (\ref{3.14.0}) we obtain
$|\nabla u|^2_{2,\Omega(\rho,\sigma)}<|u|^2_{2,\Omega(\rho,\sigma)}$. So $u\in \Theta_{\rho,\sigma}$. By Lemma \ref{l5.1.1} (i) there exists $\alpha>0$ such that $\alpha u\in\mathcal{N}_{\rho,\sigma}$.
According to the fact that $I_{\rho,\sigma}$ is weak lower semicontinuous, similar to (\ref{3.18.0}) we have
 $I_{\rho,\sigma}(\alpha u)=d_{\rho,\sigma}$.

Next we show that $I'_{\rho,\sigma}(w)=0$, where $w=\alpha u\in \mathcal{N}_{{\rho,\sigma}}$ obtained above.
Indeed, on the contrary, we assume $I'_{\rho,\sigma}(w)\neq0$. Then there exist $\delta,\nu>0$ such that
$$\|I'_{\rho,\sigma}(v)\|\geq\nu \ \ \text{for all} \ v\in E_{\rho,\sigma}\ \text{satisfying}\ \|v-w\|\leq 2\delta.$$
Since $w\neq0$, we take $L=(\inf_{\mathbb{R}^N}V)^\frac12|w|_2$ and may
assume $6\delta<L$.
Let $D=[\frac12,\frac32]$. Then as in the proof of Theorem 1.1, applying the quantitative deformation lemma and the degree theory, there will be a contradiction. Therefore, $I'_{\rho,\sigma}(w)=0$. Applying the maximum principle to (\ref{5.1}), one easily has that $w>0$. This completes the proof. \ \ \ \ $\Box$

\subsection{The existence of multiple sign-changing solutions}
For fixed $k\in\mathbb{N}$, define
$$\aligned \mathcal{N}^{\pm}_k:=\{u\in E: \ & \text{there exist}\ 0=\rho_0<\rho_1<...<\rho_k<\rho_{k+1}=+\infty\ \text{ such that}\\ &(-1)^{j}u|_{\Omega(\rho_j,\rho_{j+1})}\geq0\ \text{and} \ u|_{\Omega(\rho_j,\rho_{j+1})}\in {\mathcal{N}_{\rho_j,\rho_{j+1}}}\ \text{for } j=0,...,k\}\endaligned$$
and
$$c^{\pm}_k:=\inf_{\mathcal{N}^{\pm}_k}I.$$

We shall show that $c^{\pm}_k$ are achieved by some $u^{\pm}_k$. Since the arguments are similar, we only prove that  $c^{+}_k$ is achieved.
\begin{Remark}Similar to the characterization of (\ref{3.1.0}), we can describe $c^+_k$ as follows. Let $\mathcal{P}_k$ be the set of all functions $u\in E$ for which there
exist $0=\rho_0<\rho_1<...<\rho_k<\rho_{k+1}=+\infty$ such that $u_j:=\chi_{\Omega(\rho_j,\rho_{j+1})}u$ satisfies $0\leq (-1)^j u_j\in \Theta_{\rho_j,\rho_{j+1}}$ and $u_j\neq0$. Then
$$c^+_k=\inf_{u\in \mathcal{P}_k}\sum^{k}_{j=0}\max_{t_j>0}I(t_ju_j).$$ \end{Remark}

Now we are ready to give the proof of Theorem 1.2.

{\bf Proof of Theorem 1.2.} Firstly, we define
$$\aligned g^+(t)=\left\{ \begin{array}{lll}
g(t),\ & t\geq0,\\
-g(-t),\ & t<0,
\end{array}\right. g^-(t)=\left\{ \begin{array}{lll}
g(t),\ & t\leq0,\\
-g(-t),\ & t>0,
\end{array}\right.\endaligned$$
and $G^\pm(t)=\int^t_0g^{\pm}(s)ds$. Replacing $G$ with $G^\pm$ in the definition of $I$, we have a  new functional denoted by $I^\pm$. By Lemma \ref{l5.1.2}, the infimum
$$d^+_{\rho,\sigma}:=\inf_{\mathcal{N}_{\rho,\sigma}}I^{+}$$
is achieved. Since $|u|$ is also a minimizer, we may assume the minimizer $u$ is a positive solution of the problem (\ref{5.1}).
Similarly, we can obtain a negative minimizer for $d^-_{\rho,\sigma}:=\inf_{\mathcal{N}_{\rho,\sigma}}I^{-}$, which is a negative
solution of (\ref{5.1}).

Let $\{u_n\}$ be a minimizing sequence of $c^+_k$. We can prove that $\{u_n\}$ is bounded in $E$ by using the same arguments as in the proof of Lemma \ref{l5.1.2}. Since $u_n\in \mathcal{N}^+_k$, there exist $0=\rho^n_0<\rho^n_1<...<\rho^n_k<\rho^n_{k+1}=+\infty$ such that $(-1)^{j}{u_n}|_{\Omega(\rho^n_j,\rho^n_{j+1})}\geq0$ and ${u_n}|_{\Omega(\rho^n_j,\rho^n_{j+1})}\in \mathcal{N}_{\rho^n_j,\rho^n_{j+1}}$ for $j=0,...,k$.

{\bf Claim 1:} $\{\rho^n_{j+1}-\rho^n_j\}$ is bounded away from $0$ for each $j$.

Similar to (\ref{5.1.3}) we have
\begin{equation}\label{5.9}\varrho\leq \int_{\Omega(\rho^n_j,\rho^n_{j+1})}|f(u_n)|^{2p}\leq \int_{\Omega(\rho^n_j,\rho^n_{j+1})}|u_n|^{p},\end{equation}
where $p\in (2,2^*)$.
Moreover, using H\"{o}lder inequality and the boundedness of $\{u_n\}$ we get
\begin{equation*}\aligned \varrho\leq C\|{u_n}|_{\Omega(\rho^n_j,\rho^n_{j+1})}\|^p\bigl((\rho^n_{j+1})^N-(\rho^n_{j})^N\bigr)^{1-\frac{p}{2^*}}\leq C_1\bigl((\rho^n_{j+1})^N-(\rho^n_{j})^N\bigr)^{1-\frac{p}{2^*}}.\endaligned\end{equation*}
Then
$\{\rho^n_{j+1}-\rho^n_{j}\}$ is bounded away from $0$ for each $j$.

{\bf Claim 2:} $\{\rho^n_{k}\}$ is bounded away from $+\infty$.

If $\lim_{n\rightarrow\infty}\rho^n_{k}=+\infty$, by Strauss lemma
we have
$$u_n(x)\rightarrow0 \quad \text{as}\ |x|\rightarrow+\infty \quad \text{uniformly in }n .$$
Taking $\tilde{g}(s)=g(f(s))f'(s)+V(x)s-V(x)f(s)f'(s)$, it is easy to see that $\lim_{s\rightarrow0}\frac{\tilde{g}(s)}{s}=0$. Set $\frac{\tilde{g}(s)}{s}\Bigl|_{s=0}=0$. Thus, for fixed $\epsilon\in (0,1)$, there exists $N_0>0$ such that, for any $n\geq N_0$ there holds
$$\aligned\int_{|x|>\rho^n_{k}}(|\nabla u_n|^2+V(x)u^2_n)&=\int_{|x|>\rho^n_{k}}\tilde{g}(u_n)u_n\\
&\leq\sup_{|x|>\rho^n_{k}}\frac{|\tilde{g}(u_n)|}{|u_n|}\int_{|x|>\rho^n_{k}}u^2_n\leq \epsilon \int_{|x|>\rho^n_{k}}u^2_n.\endaligned$$
This is impossible. So we may assume $\lim_{n\rightarrow\infty}\rho^n_{k}=\rho_{k}<+\infty$.

{\bf Claim 3:} $c^+_k$ is achieved.

In fact, up to a subsequence, we may suppose that $u_n\rightharpoonup u$ in $E$, $u_n\rightarrow u$ in $L^p(\mathbb{R}^N)$ for any $p\in(2,2^*)$, and $u_n(x)\rightarrow u(x)$ a.e. in $\mathbb{R}^N$. Moreover, by Claims 1 and 2, there are $0=\rho_0<\rho_1<...<\rho_k<\rho_{k+1}=+\infty$ so that
$\rho^n_j\rightarrow\rho_j$ for $j=1,...,k$. It follows that $u_n|_{\Omega(\rho^n_j,\rho^n_{j+1})}\rightharpoonup u|_{\Omega(\rho_j,\rho_{j+1})}$ in $E$, $u_n|_{\Omega(\rho^n_j,\rho^n_{j+1})}\rightarrow u|_{\Omega(\rho_j,\rho_{j+1})}$ in $L^p(\mathbb{R}^N)$ for any $p\in(2,2^*)$, and $u_n(x)\rightarrow u(x)$ a.e. in $\Omega(\rho_j,\rho_{j+1})$. Then $(-1)^ju|_{\Omega(\rho_j,\rho_{j+1})}\geq0$. Letting $n\rightarrow\infty$ in (\ref{5.9}), we get
$$\varrho\leq C\int_{\Omega(\rho_j,\rho_{j+1})}|u|^p,$$
which implies that $u|_{\Omega(\rho_j,\rho_{j+1})}\neq0$.
 Since $u_n|_{\Omega(\rho^n_j,\rho^n_{j+1})}\in \mathcal{N}_{\rho^n_j,\rho^n_{j+1}}$ and $u_n|_{\Omega(\rho^n_j,\rho^n_{j+1})}\rightarrow u|_{\Omega(\rho_j,\rho_{j+1})}$ weakly in $E$, strongly in  $L^p(\mathbb{R}^N)$ and a.e. in $\Omega(\rho_j,\rho_{j+1})$, we deduce
 $$|\nabla u|^2_{2,\Omega(\rho_j,\rho_{j+1})}+\int_{\Omega(\rho_j,\rho_{j+1})}V(x)f(u)f'(u)u\leq \int_{\Omega(\rho_j,\rho_{j+1})}g(f(u))f'(u)u.$$
 Note that $u|_{\Omega(\rho_j,\rho_{j+1})}\neq0$, we then have
 $$|\nabla u|^2_{2,\Omega(\rho_j,\rho_{j+1})}<\int_{\Omega(\rho_j,\rho_{j+1})}g(f(u))f'(u)u<|u|^2_{2,\Omega(\rho_j,\rho_{j+1})}.$$
 Hence, $u|_{\Omega(\rho_j,\rho_{j+1})}\in \Theta_{\rho_j,\rho_{j+1}}$.
 Then there exist $\alpha_j>0$ such that \begin{equation}\label{4.9.2}\alpha_ju|_{\Omega(\rho_j,\rho_{j+1})}\in \mathcal{N}_{\rho_j,\rho_{j+1}},\end{equation} for $j=0,1,...,k$. Define
\begin{equation}\label{4.9.0}u^+_k:=\Sigma^{k}_{j=0}\alpha_ju|_{\Omega(\rho_j,\rho_{j+1})}.\end{equation}
We observe that $u^+_k\in \mathcal{N}^+_k$. Below we show that $c^+_k$ is attained. The weak convergence in $E$ and strong convergence in $L^p(\mathbb{R}^N)$ ($2<p<2^*$) of ${u_n}|_{\Omega(\rho^n_j,\rho^n_{j+1})}\rightarrow u|_{\Omega(\rho_j,\rho_{j+1})}$ imply
\begin{equation}\label{5.10}c^+_k\leq I(u^+_k)=\Sigma^{k}_{j=0}I(\alpha_ju|_{\Omega(\rho_j,\rho_{j+1})})
\leq\Sigma^{k}_{j=0}\liminf_{n\rightarrow\infty}I(\alpha_j {u_n}|_{\Omega(\rho^n_j,\rho^n_{j+1})}).\end{equation}
Moreover,
$$\Sigma^{k}_{j=0}\liminf_{n\rightarrow\infty}I(\alpha_j {u_n}|_{\Omega(\rho^n_j,\rho^n_{j+1})})\leq\Sigma^{k}_{j=0}\liminf_{n\rightarrow\infty}I( {u_n}|_{\Omega(\rho^n_j,\rho^n_{j+1})})=\liminf_{n\rightarrow\infty}I(u_n)=c^+_k.$$
Thus, $I(u^+_k)=c^+_k$. Then (\ref{5.10}) implies that
$\alpha_ju|_{\Omega(\rho_j,\rho_{j+1})}$ is a minimizer of the minimization problem
$$\inf_{\mathcal{N}_{\rho_j,\rho_{j+1}}\cap P^{\epsilon_j}}I^{\epsilon_j},$$
where $\epsilon_j=+$ for $j$ even and $\epsilon_j=-$ for $j$ odd, and $P^{\pm}:=\{u\in E:\pm u\geq0\}$. By the observation at the beginning of the proof, $\alpha_j u|_{\Omega(\rho_j,\rho_{j+1})}$ is a minimizer of
 $$\inf_{\mathcal{N}_{\rho_j,\rho_{j+1}}}I^{\epsilon_j},$$
 and for $j$ even it is a positive solution and for $j$ odd a negative solution of (\ref{5.1}),
in which $\rho=\rho_j$ and $\sigma=\rho_{j+1}$. It is obvious that  $u^+_k(\rho_j)=0$. By elliptic regularity
theory $u^+_k$ is $C^1$ on $(\rho_j,\rho_{j+1})$ for any $j$. Then the strong maximum principle
implies that $u^+_k(0)>0$, $(-1)^ju^+_k(x)>0$ for $\rho_j<|x|<\rho_{j+1}$ and $j=0,1,...,k$. So $u^+_k$ has exactly $k$ nodes.

 {\bf Claim 4:} $u^+_k$ is a solution of (\ref{2.2.0}).

 To get this claim, we shall use the method given in \cite{LW}. In fact, for simplicity, we may assume $\alpha_j=1$ in (\ref{4.9.0}) for all $j$. Argue by contradiction we assume that $I'(u^+_k)\neq0$.
Then there exists $\phi\in C^\infty_0(\mathbb{R}^N)$ such that
$ \langle I'(u^+_k),\phi\rangle=-2$. So there is $\delta>0$ such that if $|s_j-1|\leq\delta$ for
$j=0,1,...,k$ and $0\leq\epsilon\leq \delta$, then the function $\sum^k_{j=0}s_ju|_{\Omega(\rho_j,\rho_{j+1})}+\epsilon\phi$ has exactly $k$ nodes
$$0<\rho_1(s,\epsilon)<...<\rho_k(s,\epsilon)<+\infty,$$
$\rho_j(s,\epsilon)$ is continuous in $(s,\epsilon)\in D\times[0,\delta]$, where $D:=\{(s_1,...,s_k)\in\mathbb{R}^k:|s_j-1|\leq\delta\}$ and
\begin{equation}\label{5.2}\Bigl\langle I'(\sum^k_{j=0}s_j u|_{\Omega(\rho_j,\rho_{j+1})}+\epsilon\phi),\phi\Bigr\rangle<-1.\end{equation}
Define $\eta\in C(D,[0,1])$ as
$$\eta(s_1,...,s_k)=\aligned
\left\{ \begin{array}{lll}
1,\ & \text{if}\ |s_j-1|\leq \delta/4\quad \text{for all}\ j,\\
0, \ & \text{if}\ |s_j-1|\geq\delta/2\quad \text{for at least one}\ j,
\end{array}\right.\endaligned$$
Set
$$Q(s)=\sum^{k}_{j=0}s_ju|_{\Omega(\rho_j,\rho_{j+1})}+\delta\eta(s)\phi, \ \text{for all}\ s\in D.$$
Then $Q\in C(D,E)$, and for each $s\in D$, $Q(s)$ has exactly $k$ nodes
$0<\rho_1(s)<...<\rho_k(s)<+\infty$ and $\rho_j(s)$ is continuous in $s\in D$. Define $H:D\rightarrow\mathbb{R}^k$ as $H(s):=(H_1(s),...,H_k(s))$, where
$$H_j(s):=\mathcal{B}\biggl(\Bigl(\sum^k_{j=0}s_iu|_{\Omega(\rho_j,\rho_{j+1})}
+\delta\eta(s)\phi\Bigr)
\Big|_{\Omega(\rho_j(s),\rho_{j+1}(s))}\biggr),\ \text{for } j=1,...,k,$$
with $\mathcal{B}(u)=\langle I'(u),u\rangle$. Then $H\in C(D, \mathbb{R}^k)$. For fixed $j$, if $|s_j-1|=\delta$ then $\eta(s)=0$ and $\rho_i(s)=\rho_i$ for all $i=1,...,k$ and therefore $H_j(s)=\mathcal{B}(s_ju|_{\Omega(\rho_j,\rho_{j+1})})$. It then follows that $H_j(s)>0$ if $s_j=1-\delta$ and $H_j(s)<0$ if $s_j=1+\delta$. Then the degree $\deg(H,int(D),0)$ is well defined and $\deg(H,int(D),0)=(-1)^k$. Thus there is $s\in int(D)$ such that $H(s)=0$. Then $Q(s)\in \mathcal{N}^+_k$. So we have
\begin{equation}\label{5.3}I(Q(s))\geq c^+_k.\end{equation}
On the other hand, from (\ref{5.2}) we infer
$$\aligned I(Q(s))&=I\Bigl(\sum^k_{j=0}s_ju|_{\Omega(\rho_j,\rho_{j+1})}\Bigr)+\int^1_0\Big\langle I'\bigl(\sum^k_{j=0}s_j u|_{\Omega(\rho_j,\rho_{j+1})}+\theta\delta\eta(s)\phi\bigr),\delta\eta(s)\phi\Big\rangle d\theta\\&\leq I\Bigl(\sum^k_{j=0}s_ju|_{\Omega(\rho_j,\rho_{j+1})}\Bigr)-\delta\eta(s).\endaligned$$
If $|s_j-1|\leq\delta/2$ for each $j$, then
\begin{equation*} I(Q(s))<I\Bigl(\sum^k_{j=0}s_ju|_{\Omega(\rho_j,\rho_{j+1})}\Bigr)
\leq\sum^k_{j=0}I(u|_{\Omega(\rho_j,\rho_{j+1})})=c^+_k,
\end{equation*}
which contradicts (\ref{5.3}). If $|s_j-1|>\delta/2$ for at least one $j$, then \begin{equation*} I(Q(s))\leq I\Bigl(\sum^k_{j=0}s_ju|_{\Omega(\rho_j,\rho_{j+1})}\Bigr)
<\sum^k_{j=0}I(u|_{\Omega(\rho_j,\rho_{j+1})})=c^+_k,
\end{equation*} so there is still a contradiction. Hence, $I'(u_k)=0$. This ends the proof.\ \ \ \ $\Box$

\begin{Remark}\label{r4.2} Similar to the proof of Lemma \ref{l5.1.2}, the least energy $d$ given in (\ref{2.5}) is achieved by a positive function which is a ground state of the equation (\ref{2.2.0}). In the proof of Theorem 1.2, since
$\alpha_ju|_{\Omega(\rho_j,\rho_{j+1})}\in \mathcal{N}_{\rho_j,\rho_{j+1}}$ in (\ref{4.9.2}), and we understand that $\mathcal{N}_{\rho_j,\rho_{j+1}}\subset\mathcal{N}$ by defining $u(x)=0$ for $x\not\in \Omega(\rho_j,\rho_{j+1})$, we then infer $I(\alpha_ju|_{\Omega(\rho_j,\rho_{j+1})})\geq d$ and so $u^+_k$ in (\ref{4.9.0}) satisfies
$$I(u^+_k)=\sum^k_{j=0}I(\alpha_ju|_{\Omega(\rho_j,\rho_{j+1})})\geq (k+1)d.$$
\end{Remark}

\section{Proof of Theorem 1.3}
In this section, we will prove Theorem 1.3 and always assume (V$_1$), (V$_3$), (V$_4$) and ($g_3$)-($g_6$)  are satisfied.
\begin{lemma}\label{l4.1} $\mathcal{M}\neq\emptyset$.\end{lemma}
{\bf Proof}: For $i\in\{1,2\}$, let $v_i\in H^1_{0,rad}(\Omega_i)$  be the positive function satisfying
$$-\Delta v_i=\mu_i v_i\ \ \text{in}\ \Omega_i,$$
where $\mu_i$ are given in (g$_6$). Setting $\tilde{v}_i(x)=(-1)^{i-1}v_i(x)$ for $x\in \Omega_i$ and $\tilde{v}_i(x)=0$ for $x\not\in \Omega_i$. Since $\Omega_i$ has smooth boundary, we have $\tilde{v}_i\in E$ and $supp\tilde{v}_1\cap supp\tilde{v}_2=\emptyset$.

We next show that there exist $s_0,t_0>0$ such that
\begin{equation*}\langle I'(s_0\tilde{v}_1+t_0\tilde{v}_2),s_0\tilde{v}_1\rangle=\langle I'(s_0\tilde{v}_1+t_0\tilde{v}_2),t_0\tilde{v}_2\rangle=0,\end{equation*}
which implies that $s_0\tilde{v}_1+t_0\tilde{v}_2\in \mathcal{M}$ and our result follows. Indeed, for any $s,t>0$, since $\langle I'(s\tilde{v}_1+t\tilde{v}_2),s\tilde{v}_1\rangle=\langle I'(s\tilde{v}_1),s\tilde{v}_1\rangle$, we set
$\psi_1(s):=\langle I'(s\tilde{v}_1),s\tilde{v}_1\rangle$. Then
\begin{equation}\label{4.1}\aligned \psi_1(s)&=s^2|\nabla \tilde{v}_1|^2_2+\int_{\Omega_1}V(x)f(s\tilde{v}_1)f'(s\tilde{v}_1)s\tilde{v}_1-
\int_{\mathbb{R}^N}g(f(s\tilde{v}_1))f'(s\tilde{v}_1)s\tilde{v}_1\\
&=s^2\Bigl(|\nabla {v}_1|^2_{2,\Omega_1}-\int_{\Omega_1}\frac{g(f(s{v}_1))f'(s{v}_1)s{v}_1}{s^2}\Bigr).\endaligned\end{equation}
In view of (\ref{2.3}) we obtain
$$\aligned\psi_1(s)&\geq s^2|\nabla v_1|^2_{2,\Omega_1}-\epsilon s^2|v_1|^2_{2,\Omega_1}-C_{\epsilon}s^p|v_1|^p_{p,\Omega_1}\\&\geq s^2\bigl(1-\frac{\epsilon}{\mu_1}\bigr)|\nabla v_1|^2_{2,\Omega_1}-C_{\epsilon}s^p|v_1|^p_{p,\Omega_1},\endaligned$$
for some $p\in(2,2^*)$ and any $\epsilon>0$. Since $p>2$, choosing  $\epsilon=\frac{\mu_1}{2}$ we know $\psi_1(s)>0$ for small $s>0$.
In the same way as (\ref{3.4}) we have
\begin{equation}\label{4.2}\aligned&\lim_{s\rightarrow+\infty}\int_{\Omega_1}\frac{g(f(s{v}_1))f'(s{v}_1)s{v}_1}{s^2}=l\int_{\Omega_1}{v}^2_1>\mu_1\int_{\Omega_1}{v}^2_1.\endaligned\end{equation}
By (\ref{4.1}) and (\ref{4.2}) we know
$\psi_1(s)<0$ for large $s$.
Then there exists $s_0>0$ such that $\psi_1(s_0)=0$. Similarly, there exists $t_0>0$ such that $\langle I'(t_0\tilde{v}_2),t_0\tilde{v}_2\rangle=0$.
Thus $s_0\tilde{v}_1+t_0\tilde{v}_2\in \mathcal{M}$.\ \ \ \ $\Box$

\begin{lemma}\label{l4.2}$\mathcal{M}$ is closed in $E$.\end{lemma}
{\bf Proof}: By Lemma \ref{l4.1} we know $\mathcal{M}\neq\emptyset$.
As the argument in Lemma \ref{l3.4}, there exists $r>0$ such that
\begin{equation}\label{4.3}\|u^\pm\|\geq r,\ \ \forall u\in \mathcal{M}.\end{equation}Assume that there exists $\{u_n\}\subset\mathcal{M}$  such that
$u_n\rightarrow u_0$ in $E$. Then $u^\pm_n\rightarrow u^\pm_0$ in $E$ as $n\rightarrow\infty$. From (\ref{4.3}) it follows that $u^\pm_0\neq0$. Since $\langle I'(u_n),u^\pm_n\rangle=0$ and $u_n\rightarrow u_0$ in $E$, we conclude
$\langle I'(u_0),u^{\pm}_0\rangle=0$. Then $u_0\in \mathcal{M}$ and $\mathcal{M}$ is closed in $E$.\ \ \ \ $\Box$

\begin{lemma}\label{l5.3}There exists a bounded sequence $\{u_n\}\subset\mathcal{M}$ such that
$$I(u_n)\rightarrow c,\ \ \text{and}\ I'(u_n)\rightarrow0,\ \text{as}\ n\rightarrow\infty.$$
\end{lemma}
{\bf Proof:} For any $u\in \mathcal{M}$, similar to (\ref{3.10.0}) we have
$$ I(u)=I(u)-\frac12\langle I'(u),u\rangle\geq\int_{\mathbb{R}^N}\bigl[\frac14g(f(u))f(u)-G(f(u))\bigr]\geq0.$$
Thus $I$ is bounded from below on $\mathcal{M}$. Applying Ekeland's variational principle, there exists a minimizing sequence $\{u_n\}\subset \mathcal{M}$ such that
\begin{equation*} c\leq I(u_n)\leq c+\frac1n,\end{equation*}
and
\begin{equation}\label{4.6} I(v)\leq I(u_n)-\frac1n\|v-u_n\|,\ \forall v\in \mathcal{M}.\end{equation}

{\bf Claim 1}: $\{u_n\}$ is  bounded in $E$.

In fact, set $\zeta(t)=I(tu_n)$. Note that
\begin{equation*}\zeta'(t)=t\Bigl[|\nabla u_n|^2_2+\int_{\mathbb{R}^N}V(x)\frac{f(tu_n)f'(tu_n)tu_n}{t^2}-
\int_{\mathbb{R}^N}\frac{g(f(tu_n))f'(tu_n)tu_n}{t^2}\Bigr].
\end{equation*}
In view of  (\ref{3.7.0}) with $u^+$ replacing by $u_n$,  we know
$\zeta'(t)>0$ when $0<t<1$ and $\zeta'(t)<0$ when $t>1$. Hence, \begin{equation}\label{4.1.0}I(u_n)=\max_{t>0}I(tu_n).\end{equation}
Then arguing as Claim 1 in the proof of Theorem 1.1, we can show that $\{u_n\}$ is bounded in $E$. Note in particular that in the proof of (\ref{3.12.0}),  $I(u_n)\geq I(sv_n)$ is deduced by Lemma \ref{l3.1}
and here by (\ref{4.1.0}).

{\bf Claim 2}: $I'(u_n)\rightarrow0$ as $n\rightarrow\infty$.

We will borrow some ideas in \cite{NW,LIUWENXIU} to get Claim 2. Indeed, for each $n\in\mathbb{N}$ and $\varphi\in E$, setting
 $$v_n=u_n+t \varphi+su^+_n+lu^-_n.$$Define the functions $h^{\pm}_n:\mathbb{R}^3\rightarrow\mathbb{R}$ by
$$
h^{\pm}_n(t,s,l)=\int_{\mathbb{R}^N}\Bigl(|\nabla v^{\pm}_n|^2+V(x)f(v^{\pm}_n)f'(v^{\pm}_n)v^{\pm}_n-g(f(v^{\pm}_n))f'(v^{\pm}_n)v^{\pm}_n
\Bigr).$$
It is evident to show that $h^{\pm}_n(0,0,0)=0$ and $h^{\pm}_n(t,s,l)$ is of class $C^1$. Moreover, we have
$$\aligned&\frac{\partial h^{+}_n}{\partial s}(t,s,l)\\=&\int_{\mathbb{R}^N}\Bigl\{2\nabla v^+_n\nabla u^+_n+V(x)\bigl[(f'(v^+_n))^2u^+_nv^+_n+f(v^+_n)f''(v^+_n)u^+_nv^+_n+f(v^+_n)f'(v^+_n)u^+_n\bigr]
\\&-g'(f(v^+_n))(f'(v^+_n))^2u^+_nv^+_n-g(f(v^+_n))f''(v^+_n)u^+_nv^+_n-g(f(v^+_n))f'(v^+_n)u^+_n\Bigr\}.\endaligned$$
Then
\begin{equation*}\aligned&\frac{\partial h^{+}_n}{\partial s}(0,0,0)\\=&\int_{\mathbb{R}^N}\Bigl\{2(\nabla u^+_n)^2+V(x)\bigl[(f'(u^+_n))^2(u^+_n)^2+f(u^+_n)f''(u^+_n)(u^+_n)^2+f(u^+_n)f'(u^+_n)u^+_n\bigr]
\\&-g'(f(u^+_n))(f'(u^+_n))^2(u^+_n)^2-g(f(u^+_n))f''(u^+_n)(u^+_n)^2-g(f(u^+_n))f'(u^+_n)u^+_n\Bigr\}.
\endaligned\end{equation*}
Since $u_n\in \mathcal{M}$, we have
$\langle I'(u_n),u^+_n\rangle=0$. Then \begin{equation}\label{4.8}\langle I'(u^+_n),u^+_n\rangle=0.\end{equation}
Hence
\begin{equation*}\aligned&\frac{\partial h^{+}_n}{\partial s}(0,0,0)=\frac{\partial h^{+}_n}{\partial s}(0,0,0)-2\langle I'(u^+_n),u^+_n\rangle\\=&\int_{\mathbb{R}^N}V(x)\bigl[(f'(u^+_n))^2(u^+_n)^2+f(u^+_n)f''(u^+_n)(u^+_n)^2-f(u^+_n)f'(u^+_n)u^+_n\bigr]
\\&+\int_{\mathbb{R}^N}\Bigl(g(f(u^+_n))f'(u^+_n)u^+_n-g'(f(u^+_n))(f'(u^+_n))^2(u^+_n)^2-g(f(u^+_n))f''(u^+_n)(u^+_n)^2\Bigr).\endaligned\end{equation*}
In view of Lemma 2.1 (7) we infer
$$(f'(s))^2s+f(s)f''(s)s-f(s)f'(s)<0,\ \ 3f^2(s)(f'(s))^2s+f^3(s)f''(s)s-f^3(s)f'(s)>0,$$
for all $s>0$. Then
\begin{equation}\label{4.8.1}\frac{\partial h^{+}_n}{\partial s}(0,0,0)\leq \int_{\mathbb{R}^N}\left[3g(f(u^+_n))\frac{(f'(u^+_n))^2}{f(u^+_n)}(u^+_n)^2-g'(f(u^+_n))(f'(u^+_n))^2(u^+_n)^2\right]:=A^n_1.\end{equation}
By ($g_3$) and $g\in C^1(\mathbb{R})$ we have
\begin{equation}\label{5.9.0}g'(s)s-3g(s)>0, \ \ \forall s>0.\end{equation}
Then $A^n_1<0$.
Similarly, we can obtain
$$\aligned\frac{\partial h^{-}_n}{\partial l}(0,0,0)&=\int_{\mathbb{R}^N}\Bigl\{V(x)\bigl[(f'(u^-_n))^2(u^-_n)^2+f(u^-_n)f''(u^-_n)(u^-_n)^2-f(u^-_n)f'(u^-_n)u^-_n\bigr]
\\&+g(f(u^-_n))f'(u^-_n)u^-_n-g'(f(u^-_n))(f'(u^-_n))^2(u^-_n)^2-g(f(u^-_n))f''(u^-_n)(u^-_n)^2\Bigr\}<0.\endaligned$$  Moreover,
$$\frac{\partial h^{+}_n}{\partial l}(0,0,0)=0,\ \frac{\partial h^{-}_n}{\partial s}(0,0,0)=0.$$
Thus, applying the implicit function theorem, there exist $\delta_n>0$ and functions $s_n(t)$, $l_n(t)\in C^1((-\delta_n,\delta_n))$ such that $s_n(0)=l_n(0)=0$ and
\begin{equation}\label{4.9}h^{\pm}_n(t,s_n(t),l_n(t))=0,\ \ \forall t\in(-\delta_n,\delta_n),\end{equation}
which implies
\begin{equation}\label{4.9.1}v_n=u_n+t\varphi+s_n(t)u^+_n+l_n(t)u^-_n\in \mathcal{M} \ \text{for}\ t\in(-\delta_n,\delta_n).\end{equation}
We now claim that
\begin{equation}\label{5.11.0}\{s'_n(0)\} \ \text{and}\  \{l'_n(0)\}\ \text{ are bounded.}\end{equation}Indeed,
from (\ref{4.9}) it follows that
\begin{equation}\label{4.10}\aligned
\left\{ \begin{array}{lll}
\frac{\partial h^{+}_n(0,0,0)}{\partial t}+\frac{\partial h^{+}_n(0,0,0)}{\partial s}s'_n(0)+\frac{\partial h^{+}_n(0,0,0)}{\partial l}l'_n(0)=0,\\
\frac{\partial h^{-}_n(0,0,0)}{\partial t}+\frac{\partial h^{-}_n(0,0,0)}{\partial s}s'_n(0)+\frac{\partial h^{-}_n(0,0,0)}{\partial l}l'_n(0)=0,
\end{array}\right.\endaligned
\end{equation}
where
$$\aligned&\frac{\partial h^{\pm}_n}{\partial t}(0,0,0)\\=&\int_{\mathbb{R}^N}\Bigl\{2\nabla u^{\pm}_n\nabla \varphi^{\pm}+V(x)\bigl[(f'(u^{\pm}_n))^2u^{\pm}_n\varphi^{\pm}+f(u^{\pm}_n)f''(u^{\pm}_n)u^{\pm}_n\varphi^{\pm}
+f(u^{\pm}_n)f'(u^{\pm}_n)\varphi^{\pm}\bigr]
\\&-g'(f(u^{\pm}_n))(f'(u^{\pm}_n))^2u^{\pm}_n\varphi^\pm-g(f(u^\pm_n))f''(u^\pm_n)u^\pm_n\varphi^\pm
-g(f(u^\pm_n))f'(u^\pm_n)\varphi^\pm\Bigr\}.
\endaligned$$
In view of (\ref{4.10}) we get
\begin{equation}\label{4.11}s'_n(0)=\frac{1}{det D}\biggl(\frac{\partial h^-_n}{\partial t}\frac{\partial h^+_n}{\partial l}-\frac{\partial h^+_n}{\partial t}\frac{\partial h^-_n}{\partial l}\biggr)(0,0,0),\end{equation}
where the matrix
\begin{equation*}\aligned
D=\left[ \begin{array}{lll}
\frac{\partial h^{+}_n(0,0,0)}{\partial s},&\ \frac{\partial h^{+}_n(0,0,0)}{\partial l}\\
\frac{\partial h^{-}_n(0,0,0)}{\partial s},&\ \frac{\partial h^{-}_n(0,0,0)}{\partial l}
\end{array}\right].\endaligned
\end{equation*}
By (g$_5$), some computations and the boundedness of $\{u_n\}$ in $E$, we obtain
$$\biggl|\Bigl(\frac{\partial h^-_n}{\partial t}\frac{\partial h^+_n}{\partial l}-\frac{\partial h^+_n}{\partial t}\frac{\partial h^-_n}{\partial l}\Bigr)(0,0,0)\biggr|\leq C\|\varphi\|.$$
By the above fact and (\ref{4.11}), to get our claim (\ref{5.11.0}) it remains to prove that $det D$ is bounded away from zero. Indeed, if $u^+_n\rightarrow0$ in $L^r(\mathbb{R}^N)$ for any $r\in (2,2^*)$, from (\ref{2.3}) it follows that $$\int_{\mathbb{R}^N}g(f(u^+_n))f'(u^+_n)u^+_n\rightarrow0.$$
By (\ref{4.8}) we have
\begin{equation*}\label{4.12}|\nabla u^+_n|^2+\int_{\mathbb{R}^N}V(x)f^2(u^+_n)\rightarrow0.\end{equation*}
Using (\ref{2.3.4}) we know
\begin{equation*}\label{4.8.2}\|u^+_n\|\rightarrow0,\end{equation*} contradicting with (\ref{4.3}).
Therefore, by \cite[Lemma 2.1]{WM1} there exist $\delta>0$ and $y_n\in\mathbb{R}^N$ such that
$$\int_{B_1(y_n)}(u^+_n)^2>\delta.$$
Set $\bar{u}^+_n(x)={u}^+_n(x+y_n)$. Then we may assume that $\bar{u}^+_n\rightharpoonup \bar{u}^+\neq0$ in $E$. Note that $A^n_1<0$ in (\ref{4.8.1}), we deduce
$$|A^n_1|=\int_{\mathbb{R}^N}g'(f(u^+_n))(f'(u^+_n))^2(u^+_n)^2-3\int_{\mathbb{R}^N}g(f(u^+_n))\frac{(f'(u^+_n))^2}{f(u^+_n)}(u^+_n)^2.$$
In view of (\ref{5.9.0}), from Fatou lemma we know
$$\aligned&\liminf_{n\rightarrow\infty}\biggl\{\int_{\mathbb{R}^N}g'(f(u^+_n))(f'(u^+_n))^2(u^+_n)^2-3\int_{\mathbb{R}^N}g(f(u^+_n))\frac{(f'(u^+_n))^2}{f(u^+_n)}(u^+_n)^2\biggr\}
\\=&\liminf_{n\rightarrow\infty}\biggl\{\int_{\mathbb{R}^N}g'(f(\bar{u}^+_n))(f'(\bar{u}^+_n))^2(\bar{u}^+_n)^2-3\int_{\mathbb{R}^N}g(f(\bar{u}^+_n))\frac{(f'(\bar{u}^+_n))^2}{f(\bar{u}^+_n)}(\bar{u}^+_n)^2\biggr\}
\\ \geq&\int_{\mathbb{R}^N}g'(f(\bar{u}^+))f'(\bar{u}^+)(\bar{u}^+)^2-3\int_{\mathbb{R}^N}g(f(\bar{u}^+))\frac{(f'(\bar{u}^+))^2}{f(\bar{u}^+)}(\bar{u}^+)^2:=A_2>0.\endaligned$$
Then by (\ref{4.8.1}) we get
$\limsup_{n\rightarrow\infty}\frac{\partial h^{+}_n}{\partial s}(0,0,0)\leq -A_2<0$.
Similarly, there exists $A_3>0$ independent on $n$ such that
$$\limsup_{n\rightarrow\infty}\frac{\partial h^{-}_n}{\partial l}(0,0,0)\leq -A_3<0.$$
So $\liminf_{n\rightarrow\infty}det D\geq A_2A_3>0.$ Then $|s'_n(0)|\leq C\|\varphi\|$ for any $\varphi\in E$. Consequently, $\{s'_n(0)\}$ is bounded.
A similar argument can be used  to show that $\{l'_n(0)\}$ is bounded and so (\ref{5.11.0}) yields.

From (\ref{4.6}) and (\ref{4.9.1}) we have
$$\Bigl|I\bigl(u_n+t\varphi+s_n(t)u^+_n+l_n(t)u^-_n\bigr)-I(u_n)\Bigr|\leq\frac1n\|t\varphi+s_n(t)u^+_n+l_n(t)u^-_n\|.$$
Letting $t\rightarrow 0^+$ we obtain
$$|\langle I'(u_n),\varphi\rangle|\leq\frac1n \|\varphi+s'_n(0)u^+_n+l'_n(0)u^-_n\|.$$
In view of the boundedness of $\{u_n\}$, $\{s'_n(0)\}$ and $\{l'_n(0)\}$ we conclude
$$|\langle I'(u_n),\varphi\rangle|\leq o_n(1)\|\varphi\|.$$
This ends the proof.\ \ \ \  \ $\Box$

{\bf Proof of Theorem 1.3.} By Lemma \ref{l5.3}, there exists a bounded minimizing sequence $\{u_n\}\subset\mathcal{M}$ such that
$I(u_n)\rightarrow c$ {and}\ $I'(u_n)\rightarrow0$, {as} $n\rightarrow\infty.$
By the boundedness of $\{u_n\}$, we may assume that $u_n\rightharpoonup u_0$ in $E$ and $u_n\rightarrow u_0$ in $L^r(\mathbb{R}^N)$ for any $r\in (2,2^*)$. In the same way as \cite[Lemma 2.14]{HQZ}, setting $w_n=u_n-u_0$ we know
\begin{equation}\label{4.15}\langle I'(u_n),\varphi\rangle-\langle I'(u_0),\varphi\rangle-\langle I'(w_n),\varphi\rangle=o_n(1)\|\varphi\|,\end{equation}
for any $\varphi\in C^\infty_0(\mathbb{R}^N)$. It is easy to see that $I'$ is weakly sequentially continuous, and then $I'(u_0)=0$. Note that $I'(u_n)\rightarrow0$. So (\ref{4.15}) implies
$I'(w_n)\rightarrow0$. Then
$$\aligned o_n(1)&=\Bigl\langle I'(w_n),\frac{f(w_n)}{f'(w_n)}\Bigr\rangle\\&=\int_{\mathbb{R}^N}|\nabla w_n|^2[1+2f^2(w_n)(f'(w_n))^2]+\int_{\mathbb{R}^N}V(x)f^2(w_n)-\int_{\mathbb{R}^N}g(f(w_n))f(w_n).\endaligned$$
Since $w_n\rightarrow0$ in $L^r(\mathbb{R}^N)$ with $2<r<2^*$, from (\ref{2.3}) we get $\int_{\mathbb{R}^N}g(f(w_n))f(w_n)\rightarrow0$.
Hence, $$|\nabla w_n|^2_2+\int_{\mathbb{R}^N}V(x)f^2(w_n)\rightarrow0.$$
In the same way as (\ref{2.3.4}) we have $\|w_n\|\rightarrow0$.  Then
$u_n\rightarrow u_0$ in $E$ and so $I(u_0)=c$. By (\ref{4.3}), we know $u^{\pm}_0\neq0$. Then $u_0\in \mathcal{M}$. Therefore, $u_0$ is a least energy sign-changing solution of problem (\ref{2.2.0}). Similar to the argument in the proof of Theorem 1.1, we can further show that $u_0$ changes sign exactly once. This ends the proof.\ \ \ \ $\Box$

\end{document}